\documentclass{article}
\usepackage{amssymb}
\usepackage{graphicx}
\usepackage{amsmath}
\usepackage{amsfonts,amsthm}
\newtheorem{theorem}{Theorem}
\newtheorem*{thm}{Theorem 1'}

\newtheorem{conjecture}{Conjecture}
\newtheorem{corollary}{Corollary}

\newtheorem{lemma}{Lemma}

\newtheorem{proposition}{Proposition}

\def\bbr{\mathbb{R}}
\def\bbc{\mathbb{C}}
\def\bbf{\mathbb{F}}
\begin{document}
\renewcommand{\thefootnote}{}
\title{Subgroup growth of lattices in semisimple Lie groups}
\author{Alexander Lubotzky \\ Nikolay Nikolov} \date{}
\maketitle
\begin{abstract} Proving a conjecture posed in \cite{GLP}, we give
 very precise bounds for the congruence
subgroup growth of arithmetic groups. This allows us to determine
the subgroup growth of irreducible lattices of semisimple Lie
groups. In the most general case our results depend on the
Generalized Riemann Hypothesis for number fields but we can state
the following unconditional theorem:

Let $G$ be a simple  Lie group of real rank at least 2, different
than $D_4(\bbc)$, and let $\Gamma$ be any non-uniform lattice of
$G$. Let $s_n(\Gamma)$ denote the number of subgroups of index at
most $n$ in $\Gamma$. Then the limit \break $ \lim\limits_{n\to
\infty} \ \frac{\log s_n(\Gamma)}{(\log n)^2/ \log \log n}$ exists
and equals a constant $\gamma(G)$ which depends only on the Lie
type of $G$ and can be easily computed from its root system.
\end{abstract}

\section{Introduction}
\footnotetext{\emph{2000 Mathematics Subject Classification}
20H05; 22E40}

 Let $H$ be a simple real Lie group, so $H$ is the
connected part of $G(\bbr)$ for some simple algebraic group $G$.
Let $K$ be maximal compact subgroup of $H$, $X = H/K$ be the
associated symmetric space, and let $\Gamma$ be a lattice in $H$,
i.e., a discrete subgroup of finite covolume in $H$. The lattice
$\Gamma$ is said to be uniform if $H/\Gamma$ is compact and
non-uniform otherwise.  We denote by $s_n(\Gamma)$ the number of
subgroups of $\Gamma$ of index at most $n$. The study of
$s_n(\Gamma)$ for finitely generated groups $\Gamma$ has been a
focus of a lot of research in the last two decades (see \cite{LS}
and the references therein).  Our first result is a precise (and
somewhat surprising) estimate of $s_n(\Gamma)$ for higher rank
lattices.
\begin{theorem}\label{t1} Assume that $\bbr$-rank $(H) \ge 2$ and $H$  is
not locally isomorphic to $D_4(\bbc)$.  Then for every non-uniform
lattice $\Gamma$ in $H$ the limit $\lim\limits_{n\to \infty} \;
\frac{\log s_n(\Gamma)}{(\log n)^2/\log\log n} $ exists and equals
a constant $\gamma(H)$ which depends only on $H$ and not on
$\Gamma$. The number $\gamma(H)$ is an invariant which is easily
computed from the root system of $G$.
\end{theorem}
The theorem shows that different lattices in the same Lie group
have some hidden algebraic similarity; a phenomenon which also
presents itself as a corollary of Margulis super-rigidity, which
implies that $H$ can be reconstructed from each $\Gamma$.

 Every conjugacy class of subgroups of $\Gamma$ of index
$n$ has size at most $n$ (which is negligible to $s_n(\Gamma)$) and defines a unique cover of the Riemannian manifold $M =
\Gamma\backslash X$. Hence Theorem  \ref{t1} is equivalent to:
\begin{thm}\label{t1$'$} With the same assumptions on $H$ as in
Theorem \ref{t1}. Let $M$ be a finite volume non-compact manifold
covered by $X$ and let $b_n (M)$ be the number of covers of $M$ of
degree at most $n$.  Then $\lim\limits_{n\to\infty} \;
\frac{\log b_n(M)}{(\log n)^2/\log\log n}$ exists, equals $\gamma(H)$
and is independent of $M$.
\end{thm}

In spite of the geometric flavor of its statement, the proof of Theorem \ref{t1} (and
\ref{t1}') is based on a lot of number theory.
This is due to the fact that a lattice $\Gamma$ as in Theorem 1 has two
properties: \medskip

(i) \ $\Gamma$ is an arithmetic lattice by Margulis' Arithmeticity
Theorem,  and

(ii) \ $\Gamma$ has the congruence subgroup property.
\medskip

Now (i) and (ii) imply that counting finite index subgroups in
$\Gamma$ boils down to counting congruence subgroups in $\Gamma$.
In fact the main result of the current paper is the proof of the
upper bound of Conjecture \ref{conj1} below which was posed in \cite{GLP} (and one extension of the lower bounds proved there).  To
describe our results we need more terminology. \medskip

Let $G$ be a simple, simply connected, connected algebraic group
defined over a number field $k$, together with a fixed
representation $G \hookrightarrow \mathrm{GL}_{n_0}$.

Let $\mathcal{O}$ be the ring of integers of $k$. Denote by
$V_f,V_\infty$ the set of (equivalence classes of) nonarchimedean,
resp. archimedean valuations of $k$ and set $V=V_f \cup V_\infty$.
For a valuation $v\in V$, let $k_v$ denote the
completion of $k$ with respect to $v$, and similarly for $v \in V_f$ define $\mathcal{O}_v$ as the completion of $\mathcal{O}$. Let
$G_v$ be the group of $k_v$-points of $G(-)$.

Fix a finite subset $S$ of valuations of $k$ containing $V_\infty$ and
consider $\mathcal{O}_S=\{ x \in k |\quad v(x)\geq 0, \forall v
\not \in S\}$: the ring of $S$-integers of $k$. Define
$\Gamma=G(\mathcal{O}_S):=G(k)\cap
\mathrm{GL}_{n_0}(\mathcal{O}_S)$. We assume that
$G_S:=\prod_{v\in S}G_v$ is noncompact, so that $\Gamma$ is an
infinite group. \medskip

For every nonzero ideal $I$ in $\mathcal{O}_S$, let $\Gamma(I)=
\ker \left(G(\mathcal{O}_S) \rightarrow G(\mathcal{O}_S/I)
\right)$. A subgroup $\Delta$ of $\Gamma$ is called a congruence
subgroup if $\Delta$ contains $\Gamma(I)$ for some ideal $I$. Let
$C_n(\Gamma)$ be the number of congruence subgroups of $\Gamma$ of
index at most $n$. Let
\[ \alpha_+(\Gamma)= \limsup_{n \rightarrow \infty} \frac{\log C_n(\Gamma)}{(\log n)^2/ \log \log n}, \quad \textrm{and} \]
\[\alpha_-(\Gamma)= \liminf_{n \rightarrow \infty} \frac{\log C_n(\Gamma)}{(\log n)^2/ \log \log n}. \]

It was shown in \cite{GLP} that for $\Gamma=\mathrm{SL}_2(\mathbb{Z})$,
$\alpha_+(\Gamma)=\alpha_-(\Gamma)=\frac{1}{4}(3-2\sqrt{2})$. A
general conjecture was formulated there for the case where $G$ splits over $k$:

Let $R=R(G)=\frac{|\Phi_+|}{l}$, where $\Phi_+$ is the set of positive roots of the root system corresponding to $G$ and $l= \mathrm{rank}(G)$, and let
\[\gamma(G)=\frac{ (\sqrt{R(R+1)}-R)^2}{4R^2}. \] Then
\begin{conjecture}\label{conj1} $\alpha_+(\Gamma)=\alpha_-(\Gamma)=\gamma(G)$.
\end{conjecture} \bigskip

It was shown in \cite{GLP}, that assuming the Generalized Riemann
Hypothesis for Artin $L$-functions (GRH), indeed $\alpha_-(\Gamma)\geq \gamma(G)$, and
that without assuming GRH this still holds if $k/\mathbb{Q}$ is contained in an
abelian extension of $\mathbb{Q}$.

In this paper we prove the upper bound in full and extend the lower
bound result of \cite{GLP} to the non-split case. In summary:

\renewcommand{\thefootnote}{1}
\begin{theorem}\label{nonchev} Let $G$ be an absolutely simple, connected, simply connected algebraic group over a number field $k$. Let $\Phi_+$, $l$, $R(G)$ and $\gamma(G)$ are the numbers defined above for the split form of $G$. Then \medskip

\textbf{A.} $\alpha_+(\Gamma)\leq \gamma(G)$, \medskip

\textbf{B (1).} Assuming GRH we have 
\[ \alpha_-(\Gamma) \geq \gamma(G):=\frac{ (\sqrt{R(R+1)}-R)^2}{4R^2}.\]
Therefore assuming GRH it follows that $\alpha_+(\Gamma)=\alpha_-(\Gamma)=\gamma(G)$. \medskip

\textbf{B (2).} Moreover part \textbf{(1)} is unconditional provided there is a Galois field $K/\mathbb{Q}$ such that $G$ is an inner form\footnote{This term is explained in section \ref{tlb}} over $K$, and either $\mathrm{Gal}(K/\mathbb{Q})$ has an abelian subgroup of index at most 4, or deg$[K:\mathbb{Q}]<42$.
\end{theorem}

\begin{corollary} If $G$ is a Chevalley (split) group and $k=\mathbb{Q}$ then $\alpha_+(\Gamma)=\alpha_-(\Gamma)=\gamma(G)$. In particular
\[\alpha_{\pm 1}(\mathrm{SL}_d(\mathbb{Z}))= \frac{(\sqrt{d(d+2)}-d)^2}{4d^2}.\]
\end{corollary}
So Conjecture 1 is now fully proved, modulo GRH (and it is
unconditionally proved for abelian extensions $k/\mathbb{Q}$). The
case of $d=3$ of Corollary 1 was also proved independently by O.
Edhan \cite{Edhan}. 
The main content of this paper is the proof of Theorem \ref{nonchev}\textbf{A}. Part \textbf{B} is just a small improvement over \cite{GLP}.
\medskip

The extension to arbitrary $k$-simple $G$ is important when one
comes to the study of subgroup growth of lattices in a higher rank
simple Lie group $H$:

As mentioned above,  by  Margulis' Arithmeticity Theorem
(\cite{margulis}) every lattice $\Gamma$ in $H$ is arithmetic.
Moreover, a famous conjecture by Serre (\cite{serre}) asserts that
such a group $\Gamma$ has the congruence subgroup property. This
conjecture is by now proved, unless $H$ is of type $A_n$ and
$\Gamma$ is a cocompact lattice in $H$. Now, given $H$ we can
analyze the possible $G,k$ and $S$ such that $G(\mathcal{O}_S)$ is
a lattice in $H=G(\mathbb{R})^0$. The possibilities are given by
Galois cohomology and enable us to prove:
\begin{theorem}\label{t3}
Assuming  GRH and Serre's conjecture, then for every non-compact
higher rank simple Lie group $H=G(\bbr)^0$ and every lattice
$\Gamma$ in $H$ the limit \[ \lim\limits_{n\to \infty} \;
\frac{\log s_n(\Gamma)}{(\log n)^2/\log\log n}\] exists and equals
$\gamma(G)$. In particular it depends only on $H$ and not on $\Gamma$.
\end{theorem}

In fact, the proof shows that for 'most' lattices, the conclusion
of Theorem \ref{t3} holds unconditionally.  In particular this applies to
the cases treated in Theorem 1.

Theorem \ref{nonchev}\textbf{A}  was proved in \cite{GLP} in the special
case when $G=\mathrm{SL}_2$. (For general split $G$, a partial result was also obtained: $\alpha_+(\Gamma)<
C \gamma(G)$ for some absolute constant $C$.) The proof there had two parts:

(a) A reduction to an extremal problem for abelian groups (\S 5 in \cite{GLP}), and

(b) Solving this extremal problem. (Theorem 5 in \cite{GLP} restated as Theorem \ref{ab} below) \medskip

Part (a) used the explicit list of the maximal subgroups of
$\mathrm{SL}_2(\mathbb{F}_q)$. Such detailed description becomes
too long for general $G(\mathbb{F}_q)$ with the increase of the
Lie rank of $G$ and $q$.

The main new result in this work relates to part (a) and is the
following Theorem \ref{t4} (deduced in turn from its more refined
version Theorem \ref{T3} from Section \ref{reduct} below). We need some
additional notation:
\medskip

Let $X(\mathbb{F}_q)$ be a finite quasisimple group of Lie type
$X$ over the finite field $\mathbb{F}_q$ of characteristic $p>3$.
For a subgroup $H$ of $X(\mathbb{F}_q)$ denote
\[ h(H)=\frac{ \log[X(\mathbb{F}_q):H]}{\log|H^{\diamondsuit}|},\]
 where $H^{\diamondsuit}$ denotes the maximal abelian
 quotient of $H$ whose order is coprime to $p$. Set $h(H)=\infty$ if $|H^{\diamondsuit}|=1$.
\medskip

Let $\tilde{X}$ be the untwisted Lie type corresponding to $X$ (so
$X=\tilde{X}$, $^2 \tilde{X}$ or $^3 \tilde{X}$, the last case
occurring only if $\tilde{X}=D_4$). Then $\tilde{X}(-)$ is a group
scheme of a split, simple, connected algebraic group. Recall that
$R(\tilde{X})$ is the ratio of the number of positive roots of the root system of $\tilde{X}$ to its Lie rank as defined before Conjecture \ref{conj1}. Extend the
definition of $R$ to twisted Lie types by setting
$R(X)=R(\tilde{X})$.

\begin{theorem}\label{t4} Given the Lie type $X$ (twisted or untwisted) then
\[ \liminf_{q \rightarrow \infty} \ \min \{h(H)|\ H\leq X(\mathbb{F}_q) \} \geq R(X).\]
\end{theorem}

The line of the proof of Theorem \ref{t4} is the following:  We
need to minimize $h(H)$ among all subgroups of $X(\bbf_q)$. We
first show that among the parabolic subgroups the minimum (when $q
\to \infty$) is obtained for the Borel subgroup and there it is
equal to $R(X)$. (See Prop. \ref{parabolic} below).  We then show
that every  $H$ can be replaced by a parabolic subgroup ${\bf P}$
with $h({\bf P}) \le h(H) + o(1)$.  The second step itself is
divided into two stages: The case when $H$ is not contained in any
parabolic subgroup (the atomic case), and then the general case is reduced
to this case.  We stress that in this process $H$ is replaced by a
parabolic subgroup which does not necessarily contain $H$ (though
in many cases it is ``natural" and possible to choose some ${\bf
P}$ containing $H$).

 The proof of Theorem \ref{t4} does not depend on CFSG, we use
instead the work of Larsen and Pink \cite{lpink} and Liebeck, Saxl
and Seitz \cite{lss} (the latter for groups of exceptional type).

Once Theorem \ref{t4} is proved, one reduces Theorem \ref{nonchev} \textbf{A}
again to the same extremal problem on abelian groups solved in \cite{GLP}:

\begin{theorem}[Theorem 5 of \cite{GLP}]\label{ab} Let $d$ and $R\geq 1$ be fixed positive numbers. Suppose $A=C_{x_1}\times C_{x_2}\times C_{x_t}$ is an abelian group such that the orders $x_1,x_2,...,x_t$ of its cyclic factors do not repeat more than $d$ times each. Suppose that $r|A|^R\leq n$ for some positive integers $r$ and $n$. Then as $n$ tends to infinity we have
\[s_r(A)\leq n^{(\gamma +o(1))\frac{\log n}{\log \log n}},\]
where $\gamma=\frac{ (\sqrt{R(R+1)}-R)^2}{4R^2}$.
\end{theorem}

A few words about the structure of the rest of the paper: \\
In Section 2 we show how the upper bound, i.e. Theorem \ref{nonchev}\textbf{A} is
proved using Theorem \ref{T3} below, of which Theorem \ref{t4} is an
easy corollary. In Section \ref{thm8} we prove Theorem \ref{T3}. In Section \ref{tlb} we use all the previous results and
Galois cohomology to prove Theorems \ref{t1}, \ref{nonchev}\textbf{B} and \ref{t3}. We
conclude with some remarks in Section 5 relating to \cite{BGLM},
\cite{ls} and \cite{MP}. \medskip

The results of this paper are announced in \cite{GLNP}.

\section{The upper bound: reduction to Theorem \ref{T3}} \label{reduct}

\subsubsection*{Notation}

All logarithms in the paper are in base 2 unless stated otherwise. Put
\[l(n)= \frac{\log n}{\log \log n}, \quad \lambda(n)=\frac{(\log n)^2}{\log \log n}. \]

For functions $f,g$ of integral argument $n$ we write $f \sim g$ when $\frac{f(n)}{g(n)} \rightarrow 1$ as $n\rightarrow \infty$ and write $f \asymp g $ if $\log f \sim \log g$. \medskip

For a finite group $G$ we denote by $O_p(G)$ the largest normal $p$-subgroup of $G$ and $d(G)$ is the minimal size of a generating set for $G$.

The (Pr\"{u}fer) rank of $G$ is defined to be the maximal of the numbers $d(H)$ as $H$ ranges over all the subgroups of $G$. Note that this use of 'rank' is different from the $k$-rank of an algebraic group $H$, which is denoted by $rk_k(H)$.

A group $G$ is said to be a central product of its subgroups $A,B \leq G$, denoted as $G=A \circ B$, if $G=AB$ and $[A,B]=1$. \medskip

Put $\delta:=[k:\mathbb{Q}]$.

\subsubsection*{The reductions}
By our assumptions $G$ is a connected, simply connected simple algebraic group defined over $k$. Therefore there exist a finite extension $K$ of $k$  and an absolutely simple group $\bar{G}$, such that $G=\mathbf{R}_{K/k}(\bar{G})$, $G(k)=\bar{G}(K)$ and $G(\mathcal{O}_S)$ is commensurable with $\bar{G}(\bar{\mathcal{O}}_{\bar{S}})$, where $\bar{\mathcal{O}}$ is the ring of integers of $K$ and $\bar{S}$ is the set of valuations of $K$ lying above $S$. Moreover, the congruence topologies of $G(\mathcal{O}_S)$ and of $\bar{G}(\bar{\mathcal{O}}_{\bar{S}})$ are compatible. So for the purpose of counting congruence subgroups we may replace $G$ by $\bar{G}$, $K$ by $k$ and thus assume that $G$ is absolutely simple to start with.

Recall that $G$ is simply connected and $G_S$ is noncompact. Therefore by the Strong Approximation Theorem (Theorem 7.12 of \cite{pr}) the congruence subgroups of $\Gamma$ correspond to open subgroups of the cartesian product
\[ \prod_{v \in V_f \backslash S}G(\mathcal{O}_v),\] so we count subgroups of
$G(\mathcal{O}_S/I)$ for various ideals $I \vartriangleleft \mathcal{O}_S$.
\bigskip

The following is the generalization of the 'Level vs. Index' Lemma to rings of algebraic integers:

\begin{lemma}[\cite{LS}, Proposition 6.1.1] Let $H$ be a subgroup of index $n$ in $\Gamma=G(\mathcal{O}_S)$.  Then $H$ contains $\Gamma(m\mathcal{O}_S)$ for some positive integer $m\leq c_0 n$, where the constant $c_0$ depends on $G$ only.
\end{lemma}

We shall repeatedly quote results from the paper \cite{GLP}. In particular, Corollary 1.2 together with the argument in \S 1 there imply that for the upper bound it is enough to prove
\[ \limsup_{n \rightarrow \infty} \frac{\log s_n(G(\mathcal{O}_S/I_0))}{\lambda(n)}\leq \gamma(G),\]
where the ideal $I_0=m\mathcal{O}_S$ with $m\in \mathbb{N}$ satisfies  $m\leq c_0 n$.

By Corollary 6.2 of \cite{GLP} we can replace $I_0=(m)$ above with its divisor $I=\pi_1 ... \pi_t$, defined to be the product of all the different prime ideal divisors $\pi_i$ of $I_0$. Note that the norm of $I$ is at most $c'n^\delta$, where the constant $c'$ depends only on the field $k$ and the algebraic group $G$. Also $t\leq (\delta+o(1))l(n)$. \bigskip

Put \[ G_I:=\prod_i G(\mathcal{O}_S/\pi_i) \simeq G(\mathcal{O}_S/I).\]

\textbf{Remark:}
For a prime ideal $\pi$ of $\mathcal{O}_S$ belonging to a rational prime $p$ we have that $\mathcal{O}_S/\pi$ is a finite field of bounded degree: at most $\delta=[k:\mathbb{Q}]$ over $\mathbb{F}_p$. Therefore the rank of the group $G(\mathcal{O}_S/\pi)$ is bounded by a function $r=r(\dim G, k)$ of $\dim G$ and $\delta$ alone and independent of $\pi$, see Proposition 7 of Window 2 from \cite{LS}. \bigskip

Now, for a rational prime $p$ which is not coprime to $I$, (i.e $p|m$) let $M(p)$ denote the set of those ideals from $\{\pi_1,\ldots , \pi_t\}$ which divide $(p)$. Define
\[G_p:= \prod_{\pi \in M(p)}G(\mathcal{O}_S/\pi), \quad \textrm{so} \quad G_I=\prod_{p|m} G_p. \]

The strategy of the proof follows several steps in which we gradually reduce the possibilities for the subgroup $H$ of $G_I$ (each time discounting any contibutions less that $n^{o(l(n))})$:

In the first step we fix the projections $R_p$ of $H$ on each $G_p$. Then we apply the Larsen-Pink theorem to each $R_p$ which roughly says that $R_p$ resembles an algebraic subgroup. By successive reductions we deal with its unipotent part and then its semisimple part, leaving only the 'torus' (in our case just an abelian $p'$-group) as a possibility where $H$ can live. This is the point where we are in position to apply Theorem \ref{ab} and finish the proof.

While doing these reductions we need several auxilliary group-theoretic results, and in addition we have to keep track of various numerical constants (in particular the change of the index of $H$), resulting in considerable notation overload.  \bigskip

\textbf{Step 1.}

Let $R_p$ be the projection of $H \leq G_I$ on the direct factor $G_p$. We are assuming that $G$ is absolutely simple and therefore for almost all rational primes $p$ the group $G_p$ is a product of $|M(p)|\leq \delta$ quasisimple groups  $G(\mathcal{O}_S/\pi),\ \pi \in M(p)$. By the remark above it follows that the rank of $G_p$ is at most $r':=\delta r$. We deduce that there are at most $|G_p|^{r'}$ possibilities for $R_p$ in $G_p$.

Since $|G_I|=O(m^{\dim G})=O(n^{\delta \dim G})$ it follows that the number of choices for the projections $\{R_p| \ p|m \}$ is at most
\[ \prod_p |G_p|^{r'} =|G_I|^{r'}=O(n^{\delta r' \dim G }),\]
which is polynomially bounded in $n$.

Thus we can assume from now on that the set of projections  $\{R_p|\ p|m \}$ is fixed and estimate the further possibilities for $H$. \bigskip

\textbf{Step 2.}

At this stage we use the following modification of a theorem by Larsen and Pink \cite{lpink} in \cite{LP}, Corollary 3.1:
\begin{theorem}[M. Larsen, R. Pink] \label{lpink} Let $G$ be a finite subgroup of $GL_n(\mathbb{F})$, where $\mathbb{F}$ is a finite field of characteristic $p$. Then $G$ has a normal subgroup $N \geq O_p(G)$ such that:

1. $[G:N]\leq C(n)$ where $C$ depends on $n$ alone, and

2. $N/O_p(G)$ is a central product of an abelian $p'$-group $A$ and quasi-simple
groups in Lie$^*(p)$.
\end{theorem}

Here and below $Lie^*(p)$ denotes the family of finite quasisimple groups of Lie type in charactersitic $p$.

Apply Theorem \ref{lpink} to each one of the groups $R_p$: they are linear of degree at most $\delta n_0$ where $n_0$ is the degree of the linear representation of $G$. Hence there exist normal subgroups $R_p^0\geq R_p^1$ of $R_p$ such that:
\bigskip

1. $[R_p:R_p^0]\leq c$ where $c=c(n_0,\delta)$ depends on $n_0$ and $\delta$ only, and \bigskip

2. $R_p^1=O_p(R_p)$ and $R_p^0/R_p^1=A_p\circ S_p$ where $A_p$ is abelian $p'$-group and $S_p$ is quasi-semisimple of characterisitic $p$.
\bigskip

Define $R:=\prod_p R_p,\quad R^i:=\prod_p R^i_p$ for $i=0,1$, \quad $S:=\prod_p S_p$, and $A:=\prod_p A_p$. \bigskip

\textbf{Step 3.}

Consider $R_p^1$. It is a nilpotent group of nilpotency class at most $n_0$ ( By Sylow's theorem every $p$-group of $GL_{n_0}(\mathbb{F}_{p^s})$ is conjugate to a group of upper unitriangular matrices), and has rank at most $\max_p\{\textrm{rank}(G_p)\}\leq r'$. Lemma 6.1 from \cite{GLP} (also Proposition 1.3.3 in \cite{LS}) says that given $HR^1/R^1\leq R/R^1$ the number of choices for $H$ is at most
\[|R|^{3(r')^2+n_0r'}\leq (c'n^\delta)^{r'(3r'+n_0)\dim G}.\] \medskip

We are ignoring polynomial contributions to $ s_n(G(\mathcal{O}_S/I))$, therefore from now on we can assume that $H$ contains $R^1$ and count the possibilities for $\bar{H}=H/R^1$ in $\bar{R}=R/R^1$. \bigskip

\textbf{Step 4.}

The group $\bar{H}$ projects onto each factor $\bar{R}_p:=R_p/R_p^1 \geq A_p\circ S_p$ of $\bar{R}$. It follows that the nonabelian composition factors of $S_p$ counted together with their multiplicities all accur among the composition factors of $\bar{H}$. Now $\bar{R}_p/S_p$ is an abelian $p'$-group extended by a group of order at most $c$.

We \textbf{claim} that provided all primes $p$ are bigger that $c$ then $\bar{H}$ contains each $S_p$.

\textbf{Proof of claim:} Follows the proof of Lemma 4.3 in \cite{LP}: \\
Let $Z$ be the center of $S$. It is enough to show that $\bar{H}Z$ contains $S$: if so then $S=S\cap (\bar{H}Z) = (S\cap \bar{H})Z$ and therefore
$S=S \leq \bar{H}$ because $S$ is perfect. Hence we can assume that $\bar{H}$ contains $Z$ and work modulo $Z$ from now on. Note that $S/Z$ is a direct product of its factors $(ZS_p)/Z$ and they are semisimple groups over distinct fields.

Consider $\bar{H}^0:=\bar{H}\cap \bar{R}^0\leq S\circ A$. Then $\bar{H}/\bar{H}^0\simeq HR^0/R^0$ only has composition factors of order at most $c$. Therefore each simple factor of $S/Z$ (counted with its multiplicity) occurs among the composition factors of $\bar{H}^0/Z$ and hence among its derived subgroup $((\bar{H}^0)'Z)/Z\leq S/Z$.

The order of a group is the product of the orders of its composition factors.
It follows that $|((\bar{H}^0)'Z)/Z |\geq |S/Z|$ and thus $S=(\bar{H}^0)'Z$, proving the claim. $\square$ \medskip

So $\bar{H}$ contains $S$ and is thus determined by its image $\tilde{H}=\bar{H}/S$ in $\tilde{R}=\bar{R}/S$. Define $\tilde{R}^0=R^0/S$: a quotient of $A=\prod_p A_p$. \medskip

\textbf{Step 5.}

In the remaining steps  we shall reduce the problem of counting the possibilites for $\tilde{H}$ in $\tilde{R}$ to counting subgroups in certain abelian groups $E$ and $T$ (to be defined below).

The key to this reduction is the following generalization of Theorem \ref{t4}. Recall the number $R(X)$ defined in the Introduction for each Lie type $X(-)$ of simple simply connected algebraic groups over finite fields: $R(X)$ is the number of positive roots of the split form $\tilde{X}$ of $X$ divided by its rank.

\begin{theorem}\label{T3}
Let $G=X(\mathbb{F}_q)$ be a finite quasisimple group of fixed Lie type $X$ over finite field $\mathbb{F}_q$ of characteristic $p>3$. There exist a finite set $\mathcal{S}\subseteq \mathbb{Q}[x]$ of nonconstant polynomials and constants $c_1,c_2,m$ depending only on $X$ with the following property:

Suppose $H\leq G$ and $A$ is an abelian $p'$-group contained in the centre of $\bar{H}=H/O_p(H)$. Then there exist an abelian $p'$-group $T$ and a subgroup $A_0$ of $A$ such that \bigskip

(1) $A_0$ is a homomorphic image of $T$ and $[A:A_0]\leq c_1$,

(2)
\[\liminf_{q\rightarrow \infty, H\leq G} \frac{c_2+\log [G:H]}{\log |T|} \geq R(X), \quad \textrm{and} \]

(3) The group $T$ is a direct product of at most $m=m(X)$ cyclic groups, each having order $f(q)$ for some $f \in \mathcal{S}$.
\end{theorem} \medskip

For each prime ideal $\pi \in M(p)$ let $R^0_\pi$ be the projection of $R^0_p$ into the direct factor $G_\pi:= G(\mathcal{O}_S/\pi)$. Then \[ [G_p:R_p^0]\geq \prod_{\pi \in M(p)} [G_\pi:R_\pi^0]\] and $A_p$ is a subdirect product of its projections $A_\pi$ into the various $G_\pi$'s.

By our assumptions $G$ is absolutely simple. Hence for all but finitely many primes $\pi$ (which we can ignore) $G_\pi$ is a finite quasisimple group which is a form of the (split) Lie type $\tilde{X}$ of $G$. Over a finite field all the forms of $\tilde{X}$ are quasisplit and it follows that $G_\pi$ is $X(\mathcal{O}_S/\pi)$ where $X$ is a (possibly twisted) Lie type corresponding to $\tilde{X}$. For example when $G$ has type $A_n$ then $G_\pi$ is either $\mathrm{SL}_{n+1}$ or $\mathrm{SU}_{n+1}$ over finite fields. It is important to note that Theorem \ref{T3} gives the same constant $R(X)=R(\tilde{X})$ for all the forms of $G$. (In the example with $A_n$ above we have $R=(n+1)/2$.) \bigskip

Now Theorem \ref{T3} applied to $R^0_\pi \leq G_\pi$ for each $\pi \in M(p)$ gives that there is an abelian group $T_\pi$ and a subgroup $A_{\pi,0}$ of $A_\pi$ with the stated properties. In particular  $T_\pi$ maps onto $A_{\pi,0}$, and moreover
\[ [G_\pi:R_\pi^0]\geq |T_\pi|^{R(G)-o(1)}.\]
Put $A_{p,0}= \prod_{\pi \in M(p)} A_{\pi,0}$. It follows that $A_{p,0}$ is a homomorphic image of the direct product $T_p:=\prod_{\pi \in M(p)}T_\pi$ and moreover

\[[G_p:R_p^0]\geq |T_p|^{R(G)-o(1)}.\]

Define $\tilde{R}^0_p=R^0_p/S_p$. Let $E_p\leq \tilde{R}^0_p$ be the image of $A_{p,0}$ under the homomorphism $A_p\twoheadrightarrow R^0_p/S_p=\tilde{R}^0_p$. We have that $[\tilde{R}^0_p:E_p]\leq [A_p:A_{p,0}] \leq c_1^\delta$ (since $|M(p)|\leq \delta$). Also $E_p$ is a homomorphic image of $T_p$.
Let
\[T=\prod_p T_p \ \textrm{ and } \ E=\prod_p E_p. \]
Since $[G_I:R^0]=\prod_p[G_p:R^0_p]$ it now follows that \[ [G_I:R^0]\geq |T|^{R(G)-o(1)}.\]

Moreover, for any given rational prime $p$ and prime ideal $\pi$ of $\mathcal{O}_S$ dividing $p$, there are at most $\delta$ possibilities for the size of the residue field $\mathcal{O}_S/\pi$. Also there are at most $\delta$ prime ideals $\pi$ dividing $p$. We conclude that $T_p$ is a product of boundedly (by $X$ and $\delta$) many cyclic groups each having order given by a finite set of polynomials in $p$. A polynomial of degree $b>0$ cannot take the same value at more than $b$ values of its argument. Therefore there exists a number $d=d(X,\delta)$, such that the abelian group $T$ is a product of cyclic groups $C_{x_i}$ and each integer appears at most $d$ times in the sequence $\{x_i\}$.
\bigskip

\textbf{Step 6:}

We need a result which is slight generalization of Proposition 5.6 from \cite{GLP}. It allows us to pass from $\tilde{R}$ down to the abelian group $E$. We postpone its proof to Section \ref{indexpf}.

\begin{proposition} \label{index} Let $D=D_1 \times ... \times D_s$ be a direct product of finite groups, where each $D_i$ has a normal subgroup $E_i$ of index at most $C$, and $E_i$ is polycyclic of cyclic length at most $r$. Assume that the (Pr\"{u}fer) rank of each $D_i$ is at most $r$. The number of subgroups $H\leq D$ whose intersection with $E=E_1 \times ... \times E_s$ is a given subgroup $L\leq E$ is at most \[ |D|^{4r}C^{2rs^2}K, \] where $K=K(C)$ is the number of isomorphism classes of groups of order at most $C$.
\end{proposition}

Recall that the rank of each $R_p$ is at most $r'$, hence $\tilde{R}_p^0$ is abelian group of rank at most $r'$. We apply Proposition \ref{index} to $\tilde{R}=\prod_p \tilde{R}_p$ and $E=\prod_p E_p$:

Each $R_p/R_p^0$ has size at most $c$ and $[\tilde{R}^0_p:E_p]\leq c_1^{\delta}$, therefore
\[ [\tilde{R}_p:E_p]\leq [R_p:R_p^0][\tilde{R}^0_p:E_p]\leq cc_1^\delta=c_0 ,\textrm{ say.}\]
Thus, given the group $\tilde{H}\cap E$ the number of choices for $\tilde{H}$ in $\tilde{R}$ is at most \bigskip
\[|\tilde{R}|^{4r'}c_0^{2r't^2}K(c_0)^t \leq n^{4 \delta r' \dim G} c_0^{O(l(n)^2)}K^{O(l(n))}=n^{O\left(\frac{\log n}{(\log \log n)^2}\right)}=n^{o(l(n))}.\]

Since $[\tilde{R}:E] \leq c_0^t=n^{o(1)}$ it follows that $[\tilde{R}:E\cap \tilde{H}]$ and $[\tilde{R}:\tilde{H}]$ differ by at most a factor $n^{o(1)}$. So we can restrict ourselves to counting the possibilities for $\tilde{H}\cap E$. Thus without loss of generality assume that $\tilde{H} \leq E$. \bigskip

\textbf{Step 7.}

To summarize the various reductions so far: we are now counting the possibilities for $\tilde{H}\leq E$ where $E$ is a homomorphic image of $A_0=\prod_p A_{0,p}$, which is in turn an image of $T$. In turn $T=C_{x_1} \times \cdots \times C_{x_s}$ where each integer appears at most $d=d(X,\delta)$ times in the sequence $\{x_i\}$.

Let $u=[E:\tilde{H}]\leq [R^0:H]$. Then
\[ n\geq [G_I:H]=[G_I:R^0][R^0:H]\geq |T|^{R(G)-o(1)}u.\]
Hence the number of choices for $\tilde{H}$ in $E$ is at most
\[s_u(E)\leq s_u(A_0)\leq s_u(T).\]
Now we can apply Theorem \ref{ab} to the group $T$, with constant $R=R(X)$ and $d=d(X,\delta)$, giving that $s_u(T)\leq n^{(\gamma+o(1))l(n)}$. \medskip

This proves Theorem \ref{nonchev}\textbf{A} modulo Theorem \ref{ab} (proved in \cite{GLP}), Theorem \ref{T3} (proved in Section \ref{thm8}) and Proposition \ref{index}. $\square$

\subsection{Proof of Proposition \ref{index}} \label{indexpf}

We need the following
\begin{lemma} \label{cyc}
Let $A\leq B$ be groups and let $C$ and $k$ be positive integers. The number of subnormal subgroups $H$ of $B$ which contain $A$ and for which there exists a subnormal series
\[ H=H_0 \vartriangleleft H_1 \vartriangleleft ... \vartriangleleft H_k \vartriangleleft B \textrm{ with } [B:H_k]\leq C \textrm{ and } H_i/H_{i-1} \textrm{cyclic.} \]
is at most $[B:A]^k C^d K$ where $d=d(B)$ and $K$ is the number of isomorphism classes of groups of order at most $C$.
\end{lemma}

\textbf{Proof:} There are at most $K$ possibilites for the quotient group $U=B/H_k$ and then at most $|U|^d\leq C^d$ for the homomorphism $B \rightarrow U$ which determines $H_k$ as the kernel. Given $H_k$ there are at most $[H_k:A]^k$ possibilities for $H=H_0$ by Lemma 5.5 of \cite{GLP}. $\square$ \medskip

\textbf{Proof of Proposition \ref{index}:} We follow the proof of Proposition 5.6 from \cite{GLP}:

Let $F_i=D_i\times D_{i+1} \times ... \times D_s$ and let $L_i= \textrm{proj}_{F_i}L$. Denote $\tilde{L}_{i+1}=L_i \cap F_{i+1}$, so $\tilde{L}_i \leq L_i \leq F_i$.
Let $H_i= \textrm{proj}_{F_i}H$. We shall bound the number of possibiliites for the sequence $(H,H_2,...,H_s)$. 

The number of choices for $H_s \leq D_s$ is at most $|D_s|^r$ (because every subgroup of $D_i$ is generated by at most $r$ elements). Now asssume that $H_{i+1}$ is given and consider the possibilities for $H_i$. Let $X=H_i \cap F_{i+1}$, $Y=\textrm{proj}_{D_i}(H_i)$ and $Z=H_i\cap D_i$. Then $H_i$ is a subdirect product of $Y/Z$ and $H_{i+1}/X$ and is thus determined by $H_{i+1},X,Y,Z$ and isomorphism $\phi: Y/Z \rightarrow H_{i+1}/X$. 

Since rank($D_i)\leq r$ the number of choices for $Y,Z$ and $\phi$ is at most $|D_i|^r$ each. Notice that the pair of groups $H_{i+1} \geq \tilde{L}_{i+1}$ together with the group $X \geq \tilde{L}_{i+1}$ satisfies the conditions of Lemma \ref{cyc}: $H_{i+1}/X \simeq Y/Z$ and $Y/Z$ is a section of $D_i$ (so $|Y/(E_iZ \cap Y)|\leq C$ and $(E_iZ \cap Y)/Z$ is polycyclic of length $\leq r$.

Therefore the number of choices for $X$ is at most
\[[H_{i+1}:\tilde{L}_{i+1}]^r C^{d(H_{i+1})} K. \]
Now $d(H_{i+1})\leq \textrm{rank}(F_{i+1}) \leq sr$ and
\[ [H_{i+1}:\tilde{L}_{i+1}]  \leq [H_{i+1}:L_{i+1}][L_{i+1}:\tilde{L}_{i+1}] \leq C^s |D_i|,
\] because
$[H_{i+1}:L_{i+1}]\leq [H:L]\leq C^s$ and $[L_{i+1}:\tilde{L}_{i+1}]=[\textrm{proj}_{F_{i+1}}(L_i):L_i \cap F_{i+1}]\leq |D_i|$.

Thus, given $H_{i+1}$ the number of choices for $H_{i}$ is at most $|D_i|^{4r}C^{2sr}K$. Multiplying from $i=s$ to $i=1$ we obtain
\[|D|^{4r}C^{2rs^2}K^s\] as required. $\square$

\subsection{Theorem \ref{t4}}
Assuming Theorem \ref{T3}, then with the help of the Larsen-Pink result, Theorem \ref{t4} is an easy corollary:

Suppose that $H$ is subgroup of $G=X(\mathbb{F}_q)$ and let $H^0$ and $S,A\leq H^0/O_p(H)$ be the subgroups given by Larsen-Pink Theorem \ref{lpink} above. Recall that by $H^\diamondsuit$ we denote the largest abelian $p'$-quotient of $H$.

Let $L$ be the least normal subgroup of $H$ such that $H/L$ is abelian $p'$-group, then by looking at the composition factors of $L$ we see that $O_p(H)\leq L$ and then $L/O_p(H)$ must contain $S$ because the latter is a perfect group. Hence $H^{\diamondsuit}$ is a quotient of $H/S$, whence $|H^{\diamondsuit}| \leq |A|\cdot C(n)$.

Apply Theorem \ref{T3} to the group $H^0$. It follows that for some constants $c_1,c_2$ and an abelian group $T$ we have
\[\liminf_{q\rightarrow \infty, H\leq G} \frac{c_2+ \log [G:H^0]}{\log |T|} \geq R(X) \]
and $|A|\leq |T|c_1$.

Clearly $[G:H]\geq [G:H^0]/C(n)$, and together with $|H^{\diamondsuit}| \leq |T|c_1\cdot C(n)$ this easily implies the conclusion of Theorem \ref{t4}. $\square$

\section{Theorem \ref{T3}: Generalities.} \label{thm8}
Recall that $G=X(\mathbb{F}_q)$ is a finite quasisimple group of Lie type $X$ over a finite field $\mathbb{F}_q$ of characteristic $p>3$, $H$ is a subgroup of $G$ and $A$ is an abelian $p'$-group in the centre of $\bar{H}=H/O_p(H)$.
\medskip

Theorem \ref{T3} will follow from the next two Propositions:

\begin{proposition} \label{key}
In the situation of Theorem \ref{T3} there exist
constants $c_1,c_0>0$, a finite set $\mathcal{S} \leq \mathbb{Q}[x]$ of polynomials (all depending only on the Lie type $X$), an abelian $p'$-group $T$ and a parabolic subgroup $P$ of $G$, such that \bigskip

1. $T \twoheadrightarrow A_0$ for some subgroup $A_0$ of $A$ of index at most $c_1$,

2. $c_0\cdot [G:H] \geq [G:P]$ and $|T|\leq c_0\cdot |P^\diamondsuit|$, where $P^\diamondsuit$ denotes the largest abelian $p'$ -image of $P$, and

3. $T$ is a direct product of at most $m=m(X)$ cyclic groups, each having order $f(q)$ for some $f \in \mathcal{S}$.
\end{proposition}

\begin{proposition}\label{parabolic}
Let $G=X(\mathbb{F})$ be a quasisimple group of Lie type $X$ over a finite field $\mathbb{F}$ of characteristic bigger than 3. In other words $X(-)$ is an absolutely simple, connected algebraic group scheme defined over $\mathbb{F}_p$.

 Let $P=P(-)\leq X$ be a parabolic subgroup and recall the definition
\[h(H):= \frac{\log [G:H]}{\log|H^\diamondsuit|}, \quad \textrm{where } H \leq G.\] Then \[ \lim_{|\mathbb{F}|\rightarrow \infty }h(P(\mathbb{F})) \geq R(X) \] with equality if and only if $P$ is Borel subgroup of $X$.
\end{proposition}
\medskip

\textbf{Remark:} \label{isogenies} Note that given the type $X(-)$ (an
absolutely simple, connected
quiasisplit algebraic group defined over
$\mathbb{F}_p$) there are several possibilities for its fundamental
group and these give several possibilities for the finite group
$G=X(\mathbb{F})$, all of which are covers of the same finite simple
group $G/\mathrm{Z}(G)$. However a simple argument shows that once
Propositions \ref{key} and \ref{parabolic} are proved for any fixed
isogeny version of $X(-)$ they will follow for all the
others. Therefore from now on with one exception we shall assume that $X$ is simply connected and so
$G$ is the universal covering group of $G/Z$. The exception is Section \ref{classic} and the orthogonal group types
($X=B_n,D_n$ and $^2D_n$), where $X$ will be assumed to be one of
the classical groups $\Omega^{\pm}_{2n},\Omega_{2n+1}$. \medskip

Assuming the above Propositions the proof of Theorem \ref{T3} is straightforward:
Let $T$ and $P$ be the groups provided by Proposition \ref{key}. Then
 \[ \frac{c_2+\log[G:H]}{\log|T|}\geq \frac{\log[G:P]}{\log|P^\diamondsuit|}=h(P),\]
 where $c_2=2R\log c_0$. Now Proposition \ref{parabolic} gives that $\liminf \limits_{q \rightarrow \infty}h(P) \geq R(X)$ and we are done. $\square$

\subsection{Proof of Proposition \ref{parabolic}:}\label{parab}

Recall that $l$ is the untwisted Lie rank of $X$ and $\Phi_+$ is the set of positive roots. \bigskip

\textbf{Case A:} Suppose first that $X$ is untwisted Lie type.

$P(-)$ is defined by a subset of the nodes (= the fundamental roots) in the Dynkin diagram of $X$, which is a disjoint union of maximal connected subsets $C_1,C_2,...,C_n$ say, of fundamental roots. For example the following diagram defines a parabolic of $A_7(\mathbb{F})$:\\
\thicklines
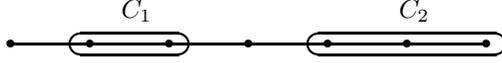
\begin{picture}(300,40)
\put (10,10){\line(1,0){180}}
\multiput(10,10)(30,0){7}{\circle*{3}}
\put(55,10){\oval(45,8)}
\put(52,20){$C_1$}
\put(160,10){\oval(75,8)}
\put(157,20){$C_2$}
\end{picture}
\bigskip

Let $E_i\subseteq \Phi_+$ be those positive roots in the span of $r\in C_i$. Then each set $E_i\cup -E_i$ is an irreducible root system with fundamental roots given by $C_i$ and Dynkin diagram which is the connected subgraph defined by $C_i$. \bigskip

Let $q=|\mathbb{F}|$. Let $L$ be the Levi factor of $P$ and let $M$ be the greatest normal subgroup of $L$ such that $L/M$ is an abelian $p'$-group. Hence $P^\diamondsuit \simeq L^\diamondsuit = L/M$. It follows that $P^\diamondsuit \simeq T/T_0$, where $T$ is a maximal split torus contained in $L$ and $T_0=M \cap T$. Since $X(-)$ is simply connected  $M$ is a direct product of its simple components and $T_0$ is also a torus. The dimension of $T_0$ is $\sum_{i=1}^n|C_i|$ and therefore
\[ [X(\mathbb{F}):P(\mathbb{F})]=q^{|\Phi_+|-\sum_{i=1}^n|E_i|}, \quad \textrm{and}\]
\[\frac{\log |P^\diamondsuit|}{\log q} \sim l-\sum_{i=1}^n|C_i| \textrm{ as } q \rightarrow \infty.\]

Notice that since $P$ is proper parabolic, the $C_i$ are proper subsets and in particular $l-\sum_{i=1}^n|C_i|>0$. It follows that
\[\lim_{q\rightarrow \infty}h(P)=\frac{|\Phi_+|-\sum_{i=1}^n|E_i|}{l-\sum_{i=1}^n|C_i|}.\]
Let $X_i$ be the split absolutely simple simply connected group having as Dynkin diagram the connected component $C_i$. Observe that $|E_i|$ is the number of positive roots of $X_i$ and $|C_i|$ is its rank. It follows that the ratio $R(X_i)$ defined in the Introduction is equal to $|E_i|/|C_i|$.

Now it is easy to check that $R(X_i)<R(X)=|\Phi_+|/l$ for every proper nonempty connected subgraph $C_i$ of the Dynkin diagram of $X$. We now use the following Lemma and obvious induction:
\begin{lemma} Suppose $a,b,c,d$ are positive real numbers such that $a>b$ and $c>d$. Suppose $\frac{a}{c}>\frac{b}{d}$. Then  $\frac{a-b}{c-d}>\frac{a}{c}$.
\end{lemma}
This shows that \bigskip
\[\frac{|\Phi_+|-\sum_{i=1}^n|E_i|}{l-\sum_{i=1}^n|C_i|}\geq R(X) \] with equality if and only if $n=0$, i.e. when $P$ is the Borel subgroup of $X(-)$.
\bigskip

\textbf{Case B:} $X$ is twisted.

We assume that char $\mathbb{F}>3$, so the corresponding untwisted type $\tilde{X}$ has Dynkin diagram with single edges and with the exception of $^3D_4$ (which can be treated similarly) $\tilde{X}$ has a symmetry $\tau$ of order 2. Also $|\mathbb{F}|=q^2$ and $\mathbb{F}$ is a quadratic extension of a field $\mathbb{F}_0$ of order $q$.
\medskip

Then $G=X(\mathbb{F})$ is the group of fixed points in $\tilde{X}(\mathbb{F})$ under the automorhism $\sigma:=\tau \phi$, where $\tau$ is the graph automorphism of $\tilde{X}(\mathbb{F})$ corresponding to the symmetry $\tau$ with the same name, and $\phi$ is the field automorphism of $\tilde{X}(\mathbb{F})$ corresponding to the automorphism $x \mapsto x^q$ of $\mathrm{Gal}(\mathbb{F}/\mathbb{F}_0)$. \medskip

The type of $G$ is $X=$$^2\tilde{X}\in \{^2A_l,^2D_l,^2E_6\}$. The root subgroups of $G$ correspond to spans $\Sigma$ of orbits of roots of $\tilde{X}$ under $\tau$ and are 1-dimensional with the exception of $\Sigma=A_2$ ocurring for $^2A_l$ with $l$ even. Here is a list of possible root subgroups:  \bigskip

\begin{tabular}{|c|l|}
Type of $\Sigma$ & Root subgroup $x_\Sigma$: \\ \hline
$A_1=\{w=w^\tau\}$ & $\{x_w(t)|\ t\in \mathbb{F}_0\}$ \\  &  \\
$A_1\times A_1=\{w,u=w^\tau\}$ & $\{x_w(t)x_u(t^q)|\ t\in \mathbb{F}\}$ \\  &  \\
$A_2=\{w,u=w^\tau,u+w\}$ & $\{x_w(t)x_u(t^q)x_{u+w}(s)|\ t,s\in \mathbb{F},$ and \\   & $t+t^q-ss^q=0$ \} \\ \hline
\end{tabular} \bigskip \\
and there is a similar parametrization for the diagonal subgroup of $G$ (see \cite{cfsg}, tables 2.4 and 2.4.7 ). \bigskip

Observe that (still excluding $^3D_4$)
\[ |x_\Sigma| \sim \sqrt{\left|\left\{ \prod_{r\in \Sigma} x_r(t_r)|\ t_r \in \mathbb{F}\right\}\right|},\]
where the right hand side is computed in $\tilde{X}(\mathbb{F})$ and the left hand side $x_\Sigma$ is a root subgroup of $G=X(\mathbb{F})$. It easily follows that
\[|G| \sim \sqrt{|\tilde{X}(\mathbb{F})|}.\]

A parabolic $P$ of $G$ is the fixed points $(\tilde{P})^\sigma$ of a parabolic $\tilde{P}$ of $\tilde{X}(\mathbb{F})$ which is defined by a $\tau$-invariant subset of the Dynkin diagram of $\tilde{X}$.

Here is an example of a parabolic of $^2A_7(-)$: \\
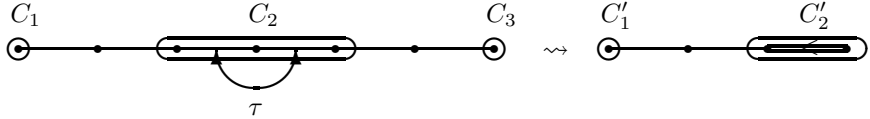
\begin{picture}(220,40)
\put (10,10){\line(1,0){180}}
\multiput(10,10)(30,0){7}{\circle*{3}}
\put(10,10){\circle{8}}
\put(7,20){$C_1$}
\put(190,10){\circle{8}}
\put(187,20){$C_3$}
\put(100,10){\oval(75,8)}
\put(97,20){$C_2$}
\put(100,10){\oval(30,30)[b]}
\put(97,-15){$\tau$}
\put(85,10){\vector(0,1){0}}
\put(115,10){\vector(0,1){0}}
\end{picture}
\begin{picture}(100,40)
\put(-15,7){$\leadsto$}
\put(10,10){\line(1,0){60}}
\put(70,9){\line(1,0){30}}
\put(70,11){\line(1,0){30}}
\multiput(10,10)(30,0){4}{\circle*{3}}
\put(10,10){\circle{8}}
\put(7,20){$C'_1$}
\put(85,10){\oval(45,8)}
\put(82,20){$C'_2$}
\put(82,8){$\mathbf{<}$}
\end{picture}
 \vspace{10mm}

From the above it easily follows that, in the notation of Case A
 \[[G:P]=\sqrt{[\tilde{X}(\mathbb{F}):\tilde{P}(\mathbb{F})]}=q^{|\Phi_+|-\sum_{i=1}^n|E_i|}, \quad \textrm{and} \]
 \[|P^\diamondsuit|\sim \sqrt{|\tilde{P}^\diamondsuit|} \sim q^{l-\sum_{i=1}^n|C_i|} \textrm{ as } q \rightarrow \infty. \]
The rest of the proof is the same as in the untwisted case. \bigskip

Finally, the case $X=$$^3D_4$ is similar to the above, with the difference that this time $|\mathbb{F}|=q^3$, $\mathbb{F}_0$ is a subfield of order $q$ and we take cube roots of the corresponding values in the untwisted group $D_4(\mathbb{F})$. $\square$

\subsection{Proof of Proposition \ref{key}.}

\subsubsection{Reduction to \emph{atomic} $H$.}

A subgroup $H$ of $G\in Lie^*(p)$ is caled \emph{p-local} if it normalizes a nontrivial $p$-subgroup of $G$.
We shall use the Borel-Tits Theorem which says that the maximal $p$-local subgroups of $G\in Lie(p)$ are parabolic:

\begin{theorem}[Borel-Tits \cite{bt}, \cite{cfsg} Theorem 3.1.1]\label{btits} Let $G\in Lie^*(p)$ be a finite quasisimple group of Lie type in characterisitc $p$ and let $R$ be a nontrivial $p$-subgroup of $G$. Then there is a parabolic subgroup $P$ of $G$, such that $R \leq O_p(P)$ and $N_G(R)\leq P$
\end{theorem}

 We shall distinguish two cases for $H$ depending on whether $H$ is $p$-local or not. We refer to the latter case as \emph{atomic}. It is the subject of sections \ref{classic} and \ref{except}. Assuming Proposition \ref{key} is proved in the atomic case, we now complete the proof in general. Thus in this section we shall
\textbf{assume that $H$ is $p$-local.} Also, since we are not interested in the explicit values of the constants $c_0,c_1$ we shall be content to define them recursively from cases of Proposition \ref{key} for type $X$ having strictly smaller Lie rank $l$. \bigskip

Now, by the Borel-Tits Theorem \ref{btits} above, we have that $H$ is contained in a proper parabolic $P'$. Choose $P'$ to be minimal parabolic containing $H$. Let $U=O_p(P')$ be the unipotent radical of $P'$, and let $L$ be its Levi factor. \bigskip

Recall that $A$ is an abelian $p'$-subgroup in the centre of $\bar{H}=H/O_p(H)$.
Thus $O_p(H)=H\cap U$ and so $\bar{H} \simeq HU/U$. We can replace $H$ by $HU$: in this way the index of $H$ in $G$ decreases, while $A$ and $\bar{H}$ stay the same (up to isomorphism). Let $H'$ be the isomorphic image of $\bar{H}$ in $L \simeq P'/U$, and identify $A$ with its isomorphic image in $H'\leq L$.

The structure of $L$ is explained in detail in Theorem 2.6.5 of \cite{cfsg}:
\begin{proposition} \label{L} Let $G=X(\mathbb{F})$ be a quasisimple group of Lie type, and let $P'$ be a parabolic subgroup of $G$ with Levi factor $L$. Define $M$ to be the largest normal subgroup of $L$ such that $L/M$ is an abelian $p'$-group (so $L/M=(P')^{\diamondsuit}$).

Then $M$ is a central product of quasisimple groups $L_1,...,L_k$ whose types correspond to connected subsets of the Dynkin diagram $X$ of $G$. When $G$ is universal (i.e. when $X(-)$ is simply connected) then each $L_i$ is universal and $M$ is in fact the direct product of the $L_i$.

Moreover there is an abelian $p'$-subgroup $T=T_L$ of $L$, such that \bigskip

1. $L_0:=MT=M \circ T$ is a central product of $T$ and $M$, \bigskip

2. $[L:L_0]\leq c_3$ for some constant $c_3$ depending only on the type $X$, and \bigskip

3. $T$ is a direct product of at most $m=m(X)$ cyclic groups whose orders are given by a finite set $\mathcal{A} \subseteq \mathbb{Q}[t]$ of nonconstant polynomials in $q$, depending on $X$ and $P$ only.
\end{proposition} \bigskip

Now, let $H_{L_0}=H'\cap L_0$ and $A_{L_0}=A\cap L_0$ so that $[A:A_{L_0}]\leq c_3$, and $A_{L_0} \leq Z(H_{L_0})$.

Put $H_M=M\cap H_{L_0}$ and $H_T=H_{L_0}\cap T$. Then $H_{L_0}/H_M$ is a quotient of $T$, so it is abelian, while $H_{L_0}/H_T$ is a quotient of $M$, so it is perfect. Therefore $H_{L_0}=H_M H_T=H_M \circ H_T$ is a central product of $H_M$ and $H_T$.

Similarly we have that $A_{L_0}= A_M \circ A_T$, where $A_M=M\cap A_{L_0}, A_T=A_{L_0} \cap T$.
\bigskip

For each direct factor $L_i$ of $M$ let $H_i$ and $A_i$ be the projections of $H_M$, resp. $A_M$  in $L_i$. Then $A_i$ is in the centre of $H_i$ and by the minimality of the parabolic $P'$ each $H_i$ is atomic in $L_i$.

The atomic case of Proposition \ref{key} applied to $A_i\leq H_i \leq L_i$ now gives that there exist constants $c(i)$, $i=1,2,...,k$, sets of nonconstant polynomials $\mathcal{S}_i \subseteq \mathbb{Q}[x]$ together with
an abelian $p'$-group $T$ and a parabolic $P_i$ of $L_i$ such that \bigskip

1. $T_i$ maps onto some subgroup $A_i(0)\leq A_i$ of index at most $c(i)$ in $A_i$, \bigskip

2. $|T_i| \leq |P_i^\diamondsuit|/c(i)$, $c(i)\cdot [L_i:H_i]\geq [L_i:P_i]$, and \bigskip

3. $T_i$ is product of boundedly many cyclic groups each having order $f(q)$ for some $f \in \mathcal{S}_i$. \bigskip

Put $T'=\prod_{i=1}^kT_i$ and $A(0)=\prod_{i=1}^k A_i(0) \leq D:=\prod_{i=1}^k A_i$. We have that $A_M$ is a subdirect product of the $A_i$, hence $A_M$ embeds in $D$, so we can identify $A_M$ as a subgroup of $D$. Put $A_M(0)=A_M \cap A(0)$. Then
\[ [A_M:A_M(0)]\leq [D:A(0)] \leq \prod_{i=1}^k c(i)=:c'_0.\]

Now $A_M(0)\leq A(0)$ and by Pontryagin duality a subgroup of a finite abelian group is also a quotient, therefore $A_M(0)$ is an image of $A(0)$, hence also an image of $T'$.

Recall that $A_{L_0}=A_M \circ A_{T}$, therefore $A_{L_0}$ is a homomorphic image of $A_M \times A_{T}$, under a map $\beta$, say.

Put $A_0=\beta (A_M(0) \times A_T)$ and $T=T'\times T_L$. Then
\[ [A:A_0]\leq [A:A_{L_0}][A_{L_0}:A_0]\leq c_3 c'_0.\]
On the other hand $A_T$ is a subgroup, hence an image of $T_L$, and therefore $A_0$ is an image of $T=T'\times T_L$. \bigskip

It is clear that $T$ satisfies condition 3 of Proposition \ref{key} for the set of polynomials $\mathcal{S}= \mathcal{A} \cup \mathcal{S}_1 \cup ... \cup \mathcal{S}_k$. It only remains to define the parabolic $P$:
\[ P:=\langle P_1,P_2,...,P_k,B \rangle ,\textrm{ where $B$ is the Borel subgroup of } G.\]
Then it is easy to see that as $q \rightarrow \infty$
\[ [G:P]\sim [G:P']\prod_{i=1}^k[L_i:P_i], \textrm{ and } |P^\diamondsuit| \sim |P'/M| \cdot \prod_{i=1}^k |P_i^\diamondsuit|. \]
Also $|T_L| \leq |P'/M| |Z(M)|$ with $|Z(M)|$ bounded by a function of $X$ alone (e.g. $2l$ where $l$ is the untwisted Lie rank of $G$).
Together with \[ [G:H] = [G:P'][L:H'] \geq [G:P'][L:H_{L_0}]/c_3 \textrm{ and} \]
\[[L:H_{L_0}] \geq \prod_{i=1}^k[L_i:H_i] \] this easily gives that condition 2 in Proposition \ref{key} is satisfied for our choice of $P$ and $T$ and appropriate constant $c_0$. \bigskip

This concludes the reduction of Proposition \ref{key} to the atomic case, i.e. when $H$ normalizes no non-trivial $p$-subgroup of $G$. This impies that every representation of $H$ over $\mathbb{F}$ is completely reducible.

\subsubsection{The atomic case I: Classical groups.} \label{classic}
By the remark on pp. \pageref{isogenies}, it is enough to prove Proposition \ref{key} for any of the isogeny versions of $X(-)$.
In this subsection we consider the case when $X$ is a classical type. Thus we may assume that $G=X(\mathbb{F})$ is one of the classical groups $\mathrm{SL}_d,\mathrm{Sp}_d, \mathrm{SU}_d,\Omega^{\pm}_{d}$ acting on its associated geometry $(V,f)$ (see Chapter 2 of \cite{KL} for the relevant definitions). Thus $V$ is a vector space of dimension $d$ over the finite field $\mathbb{F}$ with a form $f:V\times V\rightarrow \mathbb{F}$, such that one of the following holds:

(a) $f=0$, or

(b) $f$ is nondegenerate symmetric or skew-symmetric, or

(c) $f$ is nondegenerate Hermitian.\bigskip

Recall that the characteristic $p$ of $\mathbb{F}$ is assumed to be bigger than 3. In particular this avoids problems with quadratic forms in characteristic 2.

\begin{lemma} \label{irred} Suppose $U\leq V$ is an irreducible $H$-submodule. Then either $(U,f)$ is nondegenerate or else $U$ is a totally isotropic subspace for $f$.
\end{lemma}
\textbf{Proof:} The assertion is clear in case (a) when $f$ is identically 0, therefore we can assume we are in cases (b) or (c). Notice that $U^{\bot}:=\{v \in V |\ f(u,v)=0 , \forall  u \in U \}$ is an $H$-submodule of $V$ and therefore $U\cap U^{\bot}$ is a submodule of $U$. By the irreducibility of $U$ it follows that either $U\cap U^{\bot}=\{0\}$, in which case $U$ is nondegenerate, or, else $U\leq U^{\bot}$, i.e. $U$ is totally isotropic. $\square$

The parabolic subgroups of the classical groups are the stabilizers of (chains of) totally isotropic spaces. Therefore the Borel-Tits theorem implies that \bigskip

\emph{
The group $H\leq G$ is atomic if and only if $H$ stabilizes no nontrivial totally isotropic subspace of $V$.} \bigskip

In case (a) this means that $V$ is an irreducible $H$-module. In cases (b) and (c) from Lemma \ref{irred} it follows that all irreducible $H$-submodules of $V$ must be nondegenerate, and then $V$ decomposes as a direct sum
\[ V_1 \bot V_2 \bot \cdots \bot V_s \]
of pairwise orthogonal nondegenerate irreducible $H$-submodules.

Thus we are led to consider the centres of irreducible linear groups preserving a nondegenerate form $f$. In particular we have the following:

\begin{lemma}\label{atomic}
Let $H\leq \mathrm{GL}(V)$ be a finite linear group acting on a vector space $V$ of dimension $n$ over a finite field $\mathbb{F}$. Suppose that

1. $H$ is irreducible over $\mathbb{F}$, and

2. $H$ preserves a form $f: V\times V \rightarrow \mathbb{F}$ such that one of the following cases holds:

(a) $f=0$,

(b) $f$ is symmetric or skew-symmetric, bilinear and nondegenerate, or

(c) $f$ is nondegenerate Hermitian ( in which case Aut($\mathbb{F}$) is assumed to possess an involution $\sigma$). \bigskip

Then there exists a finite extension $E$ of $\mathbb{F}$ of degree $s$ say, such that $H$ is isomorphic to a group $H'\leq \mathrm{GL}(V')$, where $V'$ is an $\frac{n}{s}$-dimensional vector space over $E$ and

1. $H'$ is absolutely irreducible over $E$, i.e. $C_{\mathrm{GL}(V')}(H')=E^*$, and

2. $H'$ preserves some form $f': V'\times V' \rightarrow E$, such that

(a) $f'=0$;

(b) either  (i): $f'$ is nondegenerate bilinear symmetric or skew symmetric, or
 (ii): $f'$ is nondegenerate Hermitian and the involution $\sigma' \in Aut(E)$ fixes $\mathbb{F}$;

(c) the form $f'$ is nondegenerate Hermitian and the involution $\sigma' \in Aut(E)$ restricts to $\sigma $ on $\mathbb{F}$.
\end{lemma}
\begin{corollary}\label{cor3}
In the situation of Lemma \ref{atomic} above, let $Z$ be the centre of $H$. Then
$Z\leq E'$, where the abelian $p'$-group $E'$ is defined below for each case:

a) $E':= E^*$, a cyclic group,

(b,(i)) $E':= \{\pm 1\}$,

(b,(ii)) and (c)  $E':= \{x \in E^*|\ x^{\sigma'}=x^{-1}\}$, a cyclic group of order $\sqrt{|E|}+1$.
\end{corollary}

We delay the proof of Lemma $\ref{atomic}$ to section \ref{at}.
\bigskip

Now return to the problem.

\textbf{Case (a):} $G=\mathrm{SL}_d(q)$.

Let $E=\mathrm{End}_{\mathbb{F}H}(V)$ be the splitting field for the irreducible $H$-module $V$. Then $s=\dim_{\mathbb{F}}E$ divides the dimension $d$ of $V$. \bigskip

In this case, take \[ T=\{x\in E^*|\ \det x=(\textrm{Norm}_{E/\mathbb{F}_q}x)^{d/s}=1\}=E^* \cap \mathrm{SL}(V),\] a cyclic group of order $f_{e,s}(q)=e\frac{q^s-1}{q-1}$, where $e=(q-1,d/s)$. Again, $A$ is a subgroup, hence a quotient of $T$. Set $\mathcal{S}= \{f_{e,s}(q)=\frac{e(q^s-1)}{q-1}\ |\ e \textrm{ and } s \textrm { divide } d \}$ \bigskip

Take $A_0=A$ and define $P$ to be the stabilizer of the chain
\[\{0\}< U_1 < \cdots < U_{s-1} < V. \]
of subspaces $U_i$ with $\dim U_i=di/s$, $i=1,2,..,s-1$.

Then $\log_q|P^\diamondsuit|\sim s-1$ and $\log_q [G:P]\sim \frac{1}{2} d (d- d/s)$. \bigskip

On the other hand $H \leq \mathrm{End}_E(V) \cap \mathrm{SL}(V, \mathbb{F})$ and therefore
\[\log_q[G:H] \geq d^2-1 - (\frac{d^2}{s}-1)=d (d- \frac{d}{s}).\] \bigskip
Thus $[G:H] \geq [G:P]$ and $|T|/|P^\diamondsuit|=O(1)$ as $q \rightarrow \infty$ and we are done. \bigskip

\textbf{Case (b)}: $f$ is skew-symmetric or symmetric and $G$ is $\mathrm{Sp}_d(q)$ or $\Omega^{\pm}_d(q)$.
By Lemma \ref{atomic} the module $V$ decomposes as a sum of irreducible modules
\[ \left( V_1 \oplus \cdots \oplus V_m \right) \bigoplus \left(W_1 \oplus W_2 \oplus \cdots \oplus W_n \right), \]
where each $V_i$ has a splitting field $E_i$ and nondegenerate bilinear (symmetric or skew-symmetric) form $h_i$, say, over $E_i$ preserved by $H$. On the other hand each $W_j$ carries a nondegenerate Hermitian form $\bar{h}_j$ over its splitting field $K_j$. Let $V'_i$ (resp. $W'_j$) denote $V_i$ (resp. $W_j$) considered as vector space over $E_i$ (resp. $K_j $) together with its associated nondegenerate form $h_i$ (resp. $\bar{h}_j$).

Let $s_j=[K_j:\mathbb{F}]$, by Lemma \ref{atomic} (b,ii) the numbers $s_j$ are even and $K_j$ has an automorphism $\sigma_j$ of order 2 fixing $\mathbb{F}$.

Then $A$ acts on each irreducible $V'_i$ as $\{\pm 1\}$, and on each $W'_j$ as $\{ x \in K_j^*| \ xx^{\sigma_j}=1\}$ a cyclic group of order $f_j(q)=q^{s_j/2}+1$. Therefore it embeds in
\[ \{\pm 1 \}^m \times T, \quad \textrm{where }T:=\prod_{j=1}^n C_{f_j(q)}.\]
 We take $A_0=A\cap T$ , where $T$ is as defined above. Set $\mathcal{S}=\{f_1,f_2,...,f_d\}$ and set $c_0=2^{\dim X}$, say. We only need to define the parabolic $P$: \medskip

Observe that $H$ embeds in the direct product
\[ M:=X(V'_1) \times \cdots \times X(V'_m) \times \mathrm{U}(W'_1) \times \cdots \mathrm{U}(W'_n) \]
where $X \in \{ \mathrm{Sp} , \Omega^{\pm }\}$ as appropriate, and $\log_q |T| \sim (s_1+ ...+s_n)/2=s$, say.

Let $V_0= V_1 \oplus \cdots \oplus V_m$ and let $d_i= \dim_{\mathbb{F}_q} W'_i, \ i=1,2,...,n$. Each of the numbers $d_i$ is even. We have that
\[|\mathrm{U}(W'_i)|\sim q^{s_i/2 (d_i/s_i)^2}=q^{d_i^2/(2s_i)}.\]

Clearly $M$ is a subgroup of $M':=X(V_0) \times \mathrm{U}(W'_1) \times \cdots \mathrm{U}(W'_n)$.

Let $t=(d_1+...+d_n)/2$ and consider the chain
\[ \{0\}=U_0<U_1< \cdots U_t <V \] of $t$ totally isotropic spaces in $V$, each $U_i$ having codimension 1 in $U_{i+1}$. Let $P$ be the parabolic in $G$ which is the stabilizer of this chain. Then $|P^\diamondsuit| \sim q^t \geq q^s \sim |T|$ and we claim that $|P| \geq |M'|$:

It is easy to see that $P$ has a group isomorphic to $X(V_0)$ as a quasisimple component of its Levi factor. Moreover by its construction the unipotent part of $P$ has dimension at least equal to the number of positive roots in a root system of type $D_t$, i.e. $t(t-1)$. Hence
\[|P| \geq |X(V_0)|\cdot |P^\diamondsuit|\cdot q^{t(t-1)} \geq |X(V_0)|\cdot q^{t^2}.\]
On the other hand $|M'|=|X(V_0)| \cdot \prod_{i=1}^n|\mathrm{U}(W'_i)|$. Together with
\[\sum_{i=1}^n \log_q |\mathrm{U}(W'_i)| \leq \sum_i d_i^2/4 \leq t^2\] this justifies the claim and we are done.
\bigskip

\textbf{Case (c):}  $f$ is nondegenerate Hermitian and $G=\mathrm{SU}_d(q)$.

The $d$-dimensional $\mathbb{F}$-vector space $V$ decomposes as orthogonal sum $V=V_1\bot V_2 \bot ...\bot V_m$ of irreducible modules $V_i$. By Lemma \ref{atomic} there are finite fileds $E_i$ such that each $V_i$ is absolutely irreducible $E_iH$-module $V'_i$ which has a nondegenerate Hermitian form $h_i$ over $E_i$, preserved by $H$. Thus $H$ embeds in
\[ \mathrm{SL}(V) \cap \prod_{i=1}^m \mathrm{U}(V'_i).  \] Define $G_i=\mathrm{U}(V'_i)$ for $i=1,2,...,m$. Let $H_i, A_i$ be the projection of $H$ and $A$ into $G_i$. Then $ H \leq \mathrm{SL}(V) \cap \prod_i H_i $ and $A \leq \mathrm{SL}(V) \cap \prod_i A_i$. \bigskip

Let $d_i=\dim_{\mathbb{F}} V_i$ and let $s_i= \dim_{\mathbb{F}} E_i$. Note that by Lemma \ref{atomic} (c) all the $s_i$ must be odd. Recall that in this case $\mathbb{F}$ is a quadratic extension of a field $\mathbb{F}_0$ of order $q$. Let $\sigma'_i$ be the unique automorphism of order 2 of $E_i$ and for $i=1,2,...,m$ set
\[ T_i=\{ x\in E_i^*|\quad x^{\sigma'_i}=x^{-1}\}, \]
a cyclic group of order $q^{s_i}+1$. Since $V_i$ is absolutely irreducuble by corollary \ref{cor3} we have that $A_i \leq T_i$. Therefore
\[ A \leq M= \mathrm{SL}(V) \cap \prod_{i=1}^m T_i \]
When $i=m$ set
\[I=\{ x\in E_m^*|\quad x^{\sigma'_m}=x^{-1}\& \ \det x=(\textrm{Norm}_{E_m/\mathbb{F}}(x))^{d_m/s_m}=1 \},\] It is a cyclic group of order $f_{e,s_m}(q):=\frac{e(q^{s_m}+1)}{q+1}$, where $e=(q+1,d/s_m)$.
The index of the group $I \mathbb{F}^*=I \circ  \mathbb{F}^*$ in $T_m$ is at most $d$ and therefore by passing to a subgroup of index $\leq d$ in $A$ we may assume that
\[A \leq M':=\mathrm{SL}(V) \cap \left(T_1 \times ... \times T_{m-1} \times (I \circ \mathbb{F}^*)\right).\] Let $A_0$ be the image of $A$ under the projection $\pi: M' \rightarrow T_1 \times ... \times T_{m-1} \times I$. Then $|\ker \pi \cap A| \leq  |\mathbb{F}^*\cap I|\leq d$ and so $|A_0| \geq |A|/d$. Therefore $A_0$ is isomorphic to a subgroup of $A$ of index at most $d$ which embeds in
\[ T:=T_1 \times ... \times T_{m-1} \times I \]

As in case (a) it follows that each $H_i$ has size at most $q^{d_i^2/s_i}$, hence $\log_q |H| \leq \sum_i d_i^2/s_i -1$.
We set $T$ as above,  and
\[\mathcal{S} =\{q^{s}+1, \frac{e (q^s+1)}{q+1}\ | \quad \textrm{ $s$ and $e$ divide $d$, and $s$ is odd}\}\]

Thus $T$ is a product of at most $m \leq d$ cyclic groups whose orders are given by polynomials from  $\mathcal{S}$, and moreover $\log_q |T| \sim (\sum_i s_i)-1$.
The only thing remaining is to find an appropriate parabolic $P$ satisfying condition 2 of Proposition \ref{key}. \bigskip

Set $v= \sum_i s_i$, clearly $v \leq \sum_i d_i=d=\dim_{\mathbb{F}} V$. \bigskip

Now, consider the following chain of $[v/2]$ totally isotropic spaces in $V$:
\begin{equation}\label{chainc}
\{0\}=U_0<U_1<U_2<\cdots <U_{[v/2]},
\end{equation}
where each $U_i$ has codimension 1 in $U_{i+1}$.
Let $P$ be the parabolic stabilizing the chain (\ref{chainc}). \medskip

Now, if $v=1$ then $P=G=\mathrm{SU}_d(q)$ and we are done. Below we assume that $v\geq 2$. Then
\[\log_q|P^\diamondsuit|\sim  2[v/2] \geq v-1 \sim \log_q |T|,\]
and it is easy to see that \[ \log_q |P|\geq \frac{d(d-1)}{2}+d-1+ \frac{(d-v)(d-v-1)}{2}.\] Now use the following easy
\begin{lemma}Given positive integers $d$ and $v$, the maximum of the expression
\[ \sum_{i=1}^m \frac{d_i^2}{s_i}, \] where $s_i,d_i \in \mathbb{N}$ are subject to $s_i|d_i, d=d_1+...+d_m$ and $v=s_1+...+ s_m$ is
$(d-(v-1))^2+v-1$ and this maximum is achieved for $d_i=s_i=1$ for all $i=2,3,...,m$.
\end{lemma}
It follows that $\log_q |H| \leq (d-(v-1))^2+v-2$. Thus, in order to prove $|P|\geq |H|$ we need to check that
\[(d-(v-1))^2+v-2 \leq \frac{(d+2)(d-1)+(d-v)(d-v-1)}{2} \] which is in turn equivalent to $(v-2)d \geq v(v-3)/2$ and this inequality holds because $d\geq v \geq 2$ by our assumption. \medskip

This completes the atomic case for the classical groups.

\subsubsection{The atomic case II: \ Exceptional groups.}\label{except}
In this subsection we assume that $G=X(\mathbb{F})$ is a finite quasisimple group of exceptional type in characteristic bigger than 3, so $X\in \{E_6,E_7,E_8,\ ^2E_6,\ ^3D_4,G_2,$ $F_4\}$ ( note that $^2 D_4$ is not regarded as exceptional since it represents the orthogonal group $\Omega_8^-$).

We shall need some information on centralizers $C_G(x)$ of (non-central) semi-simple elements of $G$. The general structure theory of such centralizers is given in \cite{cfsg} Theorem 4.2.2. In our case the Lie rank $X$ of $G$ is relatively small (at most 8) so the possibilities for the components of $C_G(x)$ are quite few. In fact every such centralizer is contained either in a parabolic, or in a maximal subgroup $M$ of $G$ listed in Tables 5.1 and 5.2 of \cite{lss} (the so called groups of maximal rank).
\bigskip

Recall that in the atomic case $A$ is a subgroup of the centre of $H$. Provided $|A|$ is big enough (i.e. $|A|>K$ for some constant depending on $X$ only) then $A$ contains a semisimple element $x$ outside the center of $G$. Then $A$ lies in a maximal torus $T'$ of $G$ and $H\leq C_G(x)$. Now, in general $C_G(x)$ is either contained in a parabolic of $G$, or else it is contained in a reductive subgroup of \emph{maximal rank} of $G$, see Theorem 4.2.2 of \cite{cfsg}. However the former possibility is excluded in the atomic case.

The (maximal) subgroups of maximal rank of the exceptional Lie groups have been described by Liebeck, Saxl and Seitz, and the list can be found in tables 5.1 and 5.2 of \cite{lss}. Thus we can assume $H\leq C_G(x)\leq M$, where $M$ is one from the list in the two tables above. \bigskip

\textbf{(a) When $|M|=O(|B|)$}

Now, observe that if $|M|$ is less than a constant times the order of the Borel subgroup $B$ of $G$, then we can take the torus $T=T'$ as the required abelian group $T$ and set $A=A_0$:
We have $A \leq T$, whence $A$ is also an image of $T$ and $|T| \sim q^l$ as $q \rightarrow \infty$. Moreoever $T$ is a direct product of at most $l \leq 8$ cyclic groups each having order $f_i(q)$, where $f_i$ is from some finite set $\mathcal{S}$ of monic polynomials depending only on the type $X$ of $G$.

Clearly $H\leq M$ and if $|M|\leq c_3|B|$ for some constant $c_3$, then $ [G:H] \geq [G:B]/c_3$ and $B^\diamondsuit$ is isomorphic to the split maximal torus of $G$ hence $|B^\diamondsuit| \sim |T'|$ as $|G| \rightarrow \infty$.

Therefore $T$ and $B$ satisfy the requirements in Proposition \ref{key} for appropriate constants $c_0, c_1$.
\bigskip

\textbf{(b) When $|M|\geq |B|$}

The cases where $M$ is larger than the Borel subgroup $B$ are very few: for example the possibilities for $M$ in \cite{lss} table 5.2 are normalizers of maximal tori and have order bounded by $c\cdot q^l$ for some absolute constant $c$ (and $l$ is the Lie rank of $G$), easily giving that $|M|<|B|$.

By examining table 5.1 we list below the possibilities for the structure of those $M$ (up to conjugacy). Recall that $q=|\mathbb{F}_0|$, and let $d,e,h$ denote appropriate integers (explicitly defined in \cite{lss} but we only need that they are all bounded by an absolute constant). As usual $A.B$ denotes an extension of $B$ by $A$, and $a$ is a cyclic group of order $a$. The asymptotic $\sim$ in the last column means that as $q \rightarrow \infty$ the quantity tends to the constant specified. \bigskip

\begin{tabular}{|c|c|c|c|}
$G$ & $M$ & $\log_q |M|\sim$ & $\log_q |B|\sim $ \\ \hline
$F_4(q)$ & $d.B_4(q)$ & 36 & 28 \\  &  &  & \\
$E_6(q)$ & $h.\left(D_5(q) \times \frac{q-1}{h}\right).h$ & 46 & 42 \\  &  &  & \\
$^2E_6(q)$ & $h.\left( ^2D_5(q) \times \frac{q+1}{h}\right).h$ & 46 & 42 \\  &  &  &  \\
$E_7(q)$ & $e.\left(E_6(q) \times \frac{q-1}{e}\right).e.2$ & 79 & 70  \\
  & $e.\left( ^2E_6(q) \times \frac{q+1}{e}\right).e.2$ & 79 & 70 \\  &  &  & \\
$E_8(q)$ & $d.(A_1(q) \times E_7(q)).d$ & 136 & 128 \\
\hline
\end{tabular} \bigskip

The rest of the argument proceeds on a case by case basis: \bigskip

1. Suppose
$G=F_4(q)$ and $M=d.B_4(q)$, so $M$ is classical quasisimple group. By the argument in section \ref{classic} above applied to $H\leq M$ we can find groups $A_0$, $T$ and a parabolic $P_0$ of $M$, such that the conclusion of Proposition \ref{key} is satisfied for $H$ and $P_0$ in $M$. For example
$c[M:H] \geq [M:P_0]$, $|T|\leq c|P_0^\diamondsuit|$ etc. We use the same groups $A_0$ and $T$, and we just need to find a parabolic $P$ of $G=F_4(q)$ such that

\[ |P_0|=O(|P|), \textrm{ and } |P_0^\diamondsuit|=O(|P^\diamondsuit|) \quad \textrm{as } q \rightarrow \infty .\]

Now, there are not many possibilities for the parabolic $P_0$ in $M=B_4(q)$, and clearly if $|P_0|=O(|B|)$ then $P=B$, the Borel subgroup of $G$ will do.
It turns out that there is just one parabolic $P_0$ which fails to have order less than the Borel, and it is the largest parabolic $P_\mathrm{max}$ of $M$ which has order about $q^{29}$. However $|P_{\max}^\diamondsuit |=O(q)$ and therefore in this case we can take $P$ to be the parabolic of maximal size in $G$ (which has dimension 37 as algebraic group, and $|P^\diamondsuit| \sim q$). \bigskip

2. The rest of the cases for $M$ are even simpler:

In all of them $M$ has subgroup of 'small' (= absolutely bounded ) index which is an extension $J \rightarrow M \rightarrow C\times D$ of a direct product of two groups $C$ and $D$ by a 'small' central subgroup $J$. By going to a subgroup of small index in $H$ and then factoring $J$ we may assume that $H\leq C\times D$. Moreover $D$ is a reductive group of rank 1 (either a torus or $A_1$) and $C$ is one of the simple groups $D_5,^2D_5,E_6,^2E_6$ or $E_7$ over $\mathbb{F}$. \bigskip

Let $H_C$ and $A_C$ be the projections of $H$ and $A$ into $C$. If $A_C \not =1$, then $H_C$ is contained in $N_C(A_C)$ a subgroup of maximal rank of $C$. Therefore $[C:H_C]\geq i(C)$, where $i(C)$ is the smallest index of a subgroup of maximal rank of $C$ and $|H|\leq e\cdot |C||D|/i(C)$ for some absolute constant $e$. Now the numbers $i(C)$ for $C=E_6,^2E_6$ or $E_7$ are easy to find from table 5.1 of \cite{lss} and for $C=D_5,^2D_5$ lower bounds for $i(C)$ can be found in \cite{cooperstein}. Direct computation then shows that $|H|=o(|B|)$, so we are in the same situation as in case (a). \bigskip

Therefore we can assume that the projection of $A\leq Z(H)$ into $C$ is trivial. It follows that $A$ is a bounded extension of its intersection $A(D)=A\cap D$ with $D$, which is contained in a 1-dimensional torus $T_1$.

Thus we can select a subgroup $A_0$ of small index in $A$, which is an image of $T_1$ and for $P$ we take the parabolic of maximal size in $G$. It is certainly larger than $M$ and $P^\diamondsuit$ is one dimensional, i.e. $\log_q |P^\diamondsuit| \sim 1$ and so $P$ satisfies the conditions of Proposition \ref{key}.

This completes the proof of Proposition \ref{key} in the atomic case. $\square$
\medskip

Theorem \ref{T3} has now been proved in full.

\subsubsection{Proof of Lemma \ref{atomic}:} \label{at}
This is well-known but we were unable to find reference for it in the literature and we provide the following ad-hoc proof.

Recall that an $\mathbb{F}\, H$-module $V$ is called \emph{absolutely irreducible} if $C_{\mathrm{GL}(V)}(H)=\mathbb{F}^*$, equivalently if $V$ stays irreducible over the algebraic closure of $\mathbb{F}$.

Let $E=\mathrm{End}_{\mathbb{F}H}(V)$. By Schur's Lemma $E$ is a finite division ring and so it is a field. Say $s=[E:\mathbb{F}]$, then $V$ becomes a vector space $V'$ over $E$ of dimension $n/s$ and $G \leq \mathrm{GL}(V')$. Moreover $V'$ is an absolutely irreducible $EG$-module.

Case (a) is now finished by setting $f'=0$. For cases (b) and (c) we need to work more:

The nondegenerate form defines an antiautomorphism $A \mapsto A^*$ of $\mathrm{End}(V)$ of order 2 given by
\[ f(Au,v)=f(u,A^*v)\]
so that $A^*$ is the adjoint of $A$ with respect to $f$. It is easy to see that $E$ is stable under the adjoint map and hence it induces an automorphism $\sigma'$ of $E$ of order at most $2$. In case (b) $\sigma'$ fixes $\mathbb{F}$ while in case (c) $ \sigma'|_{\mathbb{F}}=\sigma$. Moreover as $H$ preserves the form $f$ we have that $g^*=g^{-1}$ for all $g \in H$.

Set $\epsilon=1$ unless $f$ is skew-symmetric bilinear when we set $\epsilon=-1$.
\begin{lemma} In the situation of cases (b) and (c)
there is a nondegenerate form $f':V' \times V' \rightarrow E$ and an $\mathbb{F}$-linear functional $h: E \rightarrow \mathbb{F}$ such that $h(x^{\sigma'})=\epsilon h(x)^{\sigma}$ and $f=h \circ f'$. The form $f'$ is bilinear (symmetric or skew-symmetric) if $\sigma'=1$ and is Hermitian or skew-Hermitian if $\sigma' \not=1$. More precisely we have
\begin{equation} \label{co} f'(v,u)= \epsilon f'(u,v)^{\sigma'}.\end{equation}
\end{lemma}
\textbf{Proof:} Fix $v \in V$ and define $h(x):=f(xv,v)$, it satisfies the requirements of the lemma. Now, for any pair of vectors $u,w \in V$ there is a scalar $\lambda(u,w) \in E$ such that
\[ f(xu,w)=h(\lambda(u,w) x), \quad \forall x \in E.\]
Then $\lambda(w,u)= \epsilon \lambda(u,w)^{\sigma'}$.

Let $v_1,v_2,...,v_k$ be a basis for $V=V'$ over $E$ (so $k=n/s$). Define $f'$ by
\[ f'( \sum_i \alpha_i v_i, \sum_j \beta_j v_j)= \sum_{1\leq i,j\leq k} \lambda(v_i,v_j) \alpha_i \beta_j^{\sigma'} . \]
Then (\ref{co}) is satisfied and it is easy to see that $f=h \circ f'$.
$\square$ \bigskip

We claim that $H$ preserves the form $f'$:

For a fixed $g \in H$ consider another form $f'':V' \times V' \rightarrow E$ defined by
\[f''(u,v)=f'(gu, gv)-f'(u,v). \]It is of the same type (bilinear or Hermitian) as $f'$ and
\[ h \circ f''=f(gu,gv)-f(u,v)=0.\]
Thus $f''(V',V') \subseteq \ker h<E$ giving that $f''=0$. This proves the claim. \bigskip

To finish the proof of Lemma \ref{atomic} observe that when $\epsilon=-1$ and $\sigma' \not=1$ the form $f'$ is skew-Hermitian, but we may consider instead the form $\mu f'$ where $\mu^{\sigma'}=-\mu$, and this form is Hermitian. (Recall that the characteristic of $E$ is odd and therefore such $\mu \in E$ always exists.)
$\square$
\section{The lower bound}\label{tlb} In this section we return to the notation from the Introduction, so $G$ denotes a simple, simply connected, connected algebraic group defined over a number field $k$. As explained at the beginning of Section \ref{reduct} we can further assume that $G$ is absolutely simple. Fix a linear representation of $G$, and let $\Gamma$ be an arithmetic subgroup of $G$.

The group $G$ is called $k$-quasisplit if $G$ contains a Borel subgroup defined over $k$ and $G$ is $k$-split if it contains a maximal $k$-torus which is $k$-split.

Recall that in \cite{GLP} the lower bound from Conjecture \ref{conj1} was stated and proved for split $G$. Below we show that with a little modification the same proof applies to the case when $G$ is not necessarily split. \bigskip

We shall need several basic results from Galois cohomology, which can be found in \cite{pr}, Section 2.2. Let $G_0$ be the split form of $G$ (So $G_0$ is Chevalley group of type $X$, say). Given $G_0$, then the possibilities for the $k$-isomorphism type of $G$ are parametrized by $H^1(\textrm{Gal}(\bar{k}/k), \mathrm{Aut}_k(G_0))$, the first cohomology group of the absolute Galois group Gal$(\bar{k}/k)$ with values in $\mathrm{Aut}_{\bar{k}}(G_0)$, which is usually written as $H^1(k, \mathrm{Aut}_{\bar{k}}(G_0))$.

In turn $\mathrm{Aut}_{\bar{k}}(G_0)$ is a semidirect product of $\bar{G}=G/Z(G)=G_a$, the adjoint form of $G$ by $Sym(X)$, the group of symmetries of the Dynkin diagram of $X$ preserving edge lengths:

\[ \bar{G} \longrightarrow \mathrm{Aut}_{\bar{k}}(G_0) \longrightarrow Sym(X).\]

This gives rise to the exact sequence of (noncommutative) cohomology

\[ H^1(k, \bar{G}) \longrightarrow H^1(k,\mathrm{Aut}_{\bar{k}}(G_0)) \stackrel{\alpha}{\longrightarrow} H^1(k,Sym(X)).\] \bigskip

The group Gal$(\bar{k}/k)$ acts trivially on $Sym(X)$, so that the last term is simply the conjugacy classes of homomorphisms of Gal$(\bar{k}/k)$ into $Sym(X)$. We observe that when $Sym(X)$ is non-trivial, it is usually a cyclic group of order 2, with the exception of $X=D_4$ when it is $S_3$.

The preimage of the trivial homomorphism from $H^1(k,Sym(X))$ by $\alpha$ inside $H^1(k,\mathrm{Aut}_{\bar{k}}(G_0))$ are the \emph{inner forms} of $G$, the rest are called \emph{the outer forms}. Moreover each fibre of $\alpha$ contains a  unique $k$-quasisplit representative and for inner forms this is the split form $G_0$. For example if $k'$ is a quadratic extension of $k$, the quasisplit group $\mathrm{SU}_{n+1}(k')$ is an outer $k$-form (denoted $^2 A_n$) of $X=A_n$ and the split form is $\mathrm{SL}_{n+1}(k)$. The following Proposition (to be used in section \ref{lattices}) says that we can always find an extension $E$ of very small degree over $k$, such that $G$ becomes an inner form over $E$:

\begin{proposition} \label{splitting} Let $G$ be an absolutely simple, connected, simply connected algebraic group over a number field $k$, and suppose $G$ is not a form of $D_4$. Then there exists a Galois field extension $E/k$ such that $[E:k]\leq 2$ and $G$ is an inner form over $E$.

If $G$ is form of $D_4$ then such $E$ exists with $[E:k]=1,2,3$ or $6$, the latter possibility arising only when $G$ is of type $^6D_4$.
\end{proposition}
\textbf{Proof} This is a consequence of the fact that $Sym(X)$ is a small group. Let $u \in H^1(k,\mathrm{Aut}_{\bar{k}}(G_0))$. We have to prove the existence of Galois field $E$, such that the image $\beta \circ a(u)$ in the commutative diagram below is trivial in $H^1(E,Sym(X))$: \bigskip

\[
\begin{array}{ccc}
H^1(k,\mathrm{Aut}_{\bar{k}}(G_0)) & \stackrel{\alpha}{\longrightarrow} & H^1(k,Sym(X)) \\ a \downarrow &   & b \downarrow \\
H^1(E,\mathrm{Aut}_{\bar{k}}(G_0)) & \stackrel{\beta}{\longrightarrow} & H^1(E,Sym(X))
\end{array} \] \bigskip

Now $\alpha(u) \in H^1(k,Sym(X))$ is represented by a homomorphism Gal$(\bar{k}/k) \rightarrow Sym(X)$. Let $Y\leq \textrm{Gal}(\bar{k}/k)$ be the kernel of this homomorphism and let $E$ be the fixed field of $Y$ ( so that $Y=\textrm{Gal}(\bar{k}/E)$ ). From the definition of $E$ it follows that $b\circ \alpha(u)=1= \beta \circ a(u)$ and we are done.
$\square$ \bigskip

Let $E$ be the field given by the above Proposition and
suppose $p$ is a rational prime which splits completely in $E$.
Let $\pi$ be a prime ideal of the $S$-integers $\mathcal{O}_S(E)$ of $E$ lying above $p$ and set $\pi' = \mathcal{O}_S \cap \pi$.

Then \[ \mathcal{O}_S(E)/\pi \simeq \mathcal{O}_S/\pi' \simeq \mathbb{F}_p, \]hence $G(\mathcal{O}_S/\pi')=G(\mathbb{F}_p)$ is an inner form of $G$ over $\mathbb{F}_p$.

Let the prime $p \in \mathbb{N}$ be as above. By Lang's theorem each connected algebraic group over a finite field is quasisplit, and so with the Strong Approximation Theorem we conclude that for almost all such $p$ the group $\Gamma$ maps onto the split Chevalley group $G(\mathbb{F}_p)=X(p)$ of type $X$ over $\mathbb{F}_p$. Notice that these are the same images used to prove the lower bound in \cite{GLP} in the case of \textbf{split} $G$. More precisely there it is proved the following:

\begin{theorem} \label{chev} Suppose that $G$ is a split Chevalley group, and that $k$ is contained in a Galois field $K$ over $\mathbb{Q}$.

\textbf{(i)} Assuming GRH we have 
\[ \alpha_-(\Gamma) \geq \frac{ (\sqrt{R(R+1)}-R)^2}{4R^2}.\]

\textbf{(ii)} Moreover, part \textbf{(ii)} holds unconditionally if $\mathrm{Gal}(K/\mathbb{Q})$ has an abelian subgroup of index at most 4, or if deg$[K:\mathbb{Q}]<42$.
\end{theorem} \bigskip

The proof of Theorem \ref{chev} in \cite{GLP} used only the finite images of $\Gamma$ of the form $G(\mathcal{O}_S/\pi')=G(\mathbb{F}_p)$ where $p$ is a rational prime which splits completely in $K$. Therefore the same argument proves Theorem \ref{nonchev}\textbf{B}. 

\section{Lattices in Lie groups}\label{lattices}
In this section $H$ denotes a \emph{characteristic 0 semisimple group}. By this we mean that $H=\prod_{i=1}^r G_i(K_i)$ where for each $i$, $K_i$ is a local field of characteristic 0 and $G_i$ is a connected simple algebraic group over $K_i$. The rank of $H$ is defined to be
\[\mathrm{rank}(H)= \sum_{i=1}^r \mathrm{rank}_{K_i}(G_i).\]
We assume throughout that none of the factors $G_i(K_i)$ is compact (so that
 $\mathrm{rank}_{K_i}(G_i)\geq 1$). Let $\Gamma$ be an irreducible lattice of $H$, i.e. for every infinite normal subgroup $N$ of $H$  the image of $\Gamma$ in $H/N$ is dense there.

Assume now that $\mathrm{rank}(H)\geq 2$, so by the Margulis' Arithmeticity Theorem $\Gamma$ is an $S$-arithmetic lattice. More precisely:

\begin{theorem}[16, Theorem 1]\label{margul}
There exist a number field $k$, a  connected absolutely simple
algebraic group $G$ defined over $k$, and a finite set of
valuations $S$ of $k$ containing $V_\infty$, such that $H$ is
isomorphic to $G_T= \prod_{v \in T} G_v$ for some set $T\subseteq
V$ of valuations of $k$, and moreover:

1. $\Gamma$ is the image of some $S$-arithmetic subgroup of $G$ under the embedding $G(k) \longrightarrow \prod_{v \in T} G_v$, and

2. For all $v\in S \backslash T$ the group $G_v$ is compact.
\end{theorem}

Note that the split form of $G$ is uniquely determined by the split form of the simple factors of $H$, which are necessarily of the same type. We set $\gamma (H):=\gamma(G)$, defined in the introduction for the split form of the algebraic group $G$.

Since $\Gamma$ is commensurable with $G(\mathcal{O}_S)$ the two
groups have 'roughly the same' subgroup growth. This statement
can be made precise, see Proposition 1.11.1 of \cite{LS}. Passing
to the simply connected cover of $G$ also does not affect the
asymptotics of the subgroup growth (see Proposition 1.11.2 of
\cite{LS}), and therefore we can assume that $G$ is in fact simply
connected. As $S$-rank$(G)= \mathrm{rank}(H) \geq 2$,  Serre's
conjecture,  (on the finiteness of the congruence kernel of
$G(\mathcal{O}_S)$, see \cite{serre}) gives that the congruence
subgroup growth of $G(\mathcal{O}_S)$ is (asymptotically) the same
as its subgroup growth. \bigskip

Now the results of the previous sections (which rely on the GRH at one point: Theorem \ref{nonchev}\textbf{B}) imply that

\[\lim_{n \rightarrow \infty} \frac{\log s_n(\Gamma)}{(\log n)^2/ \log \log n}= \lim_{n \rightarrow \infty} \frac{\log C_n\left(G(\mathcal{O}_S)\right)}{(\log n)^2/ \log \log n}=\gamma(G).\] \bigskip

Thus Theorem \ref{t3} is now proved modulo the validity of the
Generalized Riemann Hypothesis for number fields and Serre's
conjecture, on the finiteness of the congruence kernel. In fact we
have proved more:

\begin{theorem}\label{t11} Let $H$ be a semisimple group with rank$(H)\geq 2$. Assuming $GRH$ and Serre's conjecture, then  for every irreducible lattice $\Gamma$ of $H$ the limit
\[ \lim_{n \rightarrow \infty} \frac{ \log s_n(\Gamma)}{(\log n)^2/\log \log n}
\] exists and equals $\gamma(H)$, i.e. it depends only on $H$ and not on $\Gamma$.
\end{theorem}

\subsection{Theorem \ref{t1}}
 When $H$ is \emph{simple} and not locally isomorphic to $D_4(\mathbb{C})$ and $\Gamma$ is a non-uniform lattice ( i.e. $\Gamma \backslash H$ is noncompact ) we can remove the dependence on GRH and Serre conjecture above: \bigskip

Indeed then $T$ must consist of a single valuation, and as $\Gamma$ is non-uniform $G$ is $k$-isotropic. Therefore $G_v$ is never compact for any $v \in V$. It follows that $S=T$, in particular $k$ has only one archimedean valuation. Hence $k$ is either $\mathbb{Q}$ or an imaginary quadratic extension of $\mathbb{Q}$.

Recall that with the exception of $G=$ $^6D_4$ the extension $E$ given by Proposition \ref{splitting} has degree at most 3 over $k$. In that case the Galois closure $K$ of $E$ over $\mathbb{Q}$ is rather small: Its Galois group $\Delta:=\textrm{Gal}(K/\mathbb{Q})$ has subnormal series
\[ \Delta \vartriangleright \Delta_1 \vartriangleright \Delta_2,\]
where $[\Delta:\Delta_1]\leq 2$, $[\Delta_1: \Delta_2]\leq 3$ and $\Delta_2$ is core-free in $\Delta$.
 An easy group theoretic argument now gives that $|\Delta|=[K:\mathbb{Q}]$ divides 18 or 8 and for such fields $K$, Theorem \ref{nonchev} part \textbf{B(2)} is true unconditionally. \bigskip

When $G$ is $^6 D_4$ then $E/k$ may have degree 6 and Galois group $S_3$ and then $[K:\mathbb{Q}]$ divides 72. The only case not covered by Theorem \ref{nonchev} part \textbf{B(2)} is when the degree is exactly 72. Indeed this is the reason that we exclude $D_4(\mathbb{C})$: In this case we must have that $k$ is an imaginary quadratic extension of $\mathbb{Q}$, so $H$ is locally isomorphic to $D_4(\mathbb{C})$. If the form of $\Gamma$ comes from a form of type $^6 D_4$ we need to use the GRH. For the other lattices in $D_4(\mathbb{C})$ the result is true unconditionally.

\medskip

Finally, note that when $G$ is $k$-isotropic the truth of Serre's conjecture has been verified: see Theorem 9.17 of \cite{pr}. \medskip

Theorem \ref{t1} is now clear. $\square$

\section{Concluding remarks}

Let us relate the results of this paper with those of \cite{BGLM} on one hand and those of \cite{ls} and \cite{MP} on the other hand. \medskip

Theorem \ref{t11} above gives a very precise estimate for the subgroup growth of lattices in higher rank semisimple groups. By way of contrast, when $H$ is of rank 1 then the type of growth is in general very different: type $n^n$ instead of $n^{\log n/ \log \log n}$. (See \cite{LS} Chapter 7.2 for a detailed discussion; only partial results are known.)

Thus, if rank$(H)=1$ and $\Gamma \leq H$ is a lattice it is natural to try to study the asymptotic behaviour of $\log s_n(\Gamma)/(n\log n)$. The following result has been proved recently independently by Liebeck and Shalev and by M\"{u}ller and Puchta:
\begin{theorem}[\cite{ls}, \cite{MP}] \label{t12}
If $H= \mathrm{PSL}_2(\mathbb{R})$ and $\Gamma$ is a lattice in $H$ then
\[ \lim_{n \rightarrow \infty} \frac{ \log s_n(\Gamma)}{\log n!}= -\chi (\Gamma)\]
where $\chi(\Gamma)$ denotes the Euler characteristic of $\Gamma$.
\end{theorem}

The proof of Theorem \ref{t12} relies on the explicit known presentations of lattices in $\mathrm{PSL}_2(\mathbb{R})$ (which are Fuchsian groups). Thus one cannot expect these methods to work for the general rank 1 groups. They still may be extended to the case of rank one groups over nonarchimedean local fields. For such an $H$ every lattice is cocompact and virtually free. The group $H=\mathrm{PSL}_2(\mathbb{Q}_p)$ is an interesting first test case. For some explicit presentations of lattices there see \cite{LM}.

We should mention however that Theorem \ref{t12} in its current form is not true for general lattices in other rank one simple groups. Indeed if $H=\mathrm{PSL}_2(\mathbb{C})$ and $\Gamma$ is a cocompact subgroup of $H$ then it follows from Poincar\'{e} duality that $\chi(\Gamma)=0$. On the other hand there exist cocompact lattices in $\mathrm{PSL}_2(\mathbb{C})$ which are mapped onto a non-abelian free groups, see \cite{lu}. For such lattices clearly $\lim_{n \rightarrow \infty} \frac{ \log s_n(\Gamma)}{\log n!}$ is positive, if it exists. A similar remark applies to $\mathrm{SO}(n,1)$, when $n$ is odd. (Note that $\mathrm{PSL}_2(\mathbb{C})$ is locally isomorphic to $\mathrm{SO}(3,1)$.)

Recall that with a suitable normalization of the Haar measure on
$\mathrm{PSL}_2(\mathbb{R})$, for every lattice $\Gamma$ in
$\mathrm{PSL}_2(\mathbb{R})$ we have
$-\chi(\Gamma)=\mathrm{vol}(\mathrm{PSL}_2(\mathbb{R})/\Gamma)$.
One may speculate and suggest that for a general lattice $\Gamma$
in $G=\mathrm{PSL}_2(\mathbb{C})$ (or $G=\mathrm{SO}(n,1)$) the
limit $\lim_{n \rightarrow \infty} \frac{ \log s_n(\Gamma)}{\log
n!}$ exists and is proportional to the covolume of $\Gamma$ in
$G$. This may be a possible way to extend Theorem \ref{t12} to
more general rank 1 groups.
\medskip

It is also of interest to relate the results of the current paper
to those of \cite{BGLM}. There,  the following invariant of a
simple Lie group $H$ was studied: For $r \in \mathbb{R}_+$ denote
by $\alpha_H(r)$ the number of conjugacy classes of lattices of
$H$ of covolume at most $r$. By a result of Wang this number is
finite if $H$ is not $\mathrm{PSL}_2(\mathbb{R})$ or
$\mathrm{PSL}_2(\mathbb{C})$. It is proved in $\cite{BGLM}$ that
for $H=\mathrm{SO}(d,1),\ d \geq 4$ there exist two positive
constants $a(d)$ and $b(d)$ such that
\[ a(d) r \log r \leq \log \alpha_H(r) \leq b(d) r \log r
\] for all sufficiently large $r$. It is further conjectured there that for simple Lie groups $H$ of higher rank there exist $a(H)$ and $b(H)$ such that
\[ a(H) \frac{(\log r)^2}{\log \log r} \leq \log \alpha_H(r) \leq b(H) \frac{(\log r)^2}{\log \log r}. \]

The results of the current paper support a stronger conjecture: the limit
\[\lim_{n \rightarrow \infty} \frac{\log \alpha_H(r)}{(\log r)^2/\log \log r}\]
exists and equals $\gamma (H)$.
\subsubsection*{Acknowledgements:}  The authors thank Inna Korchagina for helpful
discussion of the subgroups of the simple groups of exceptional
type. This research was carried out while the second author held a
Golda-Meir Postdoctoral Fellowship at the Hebrew University of
Jerusalem. The support of the Israel Science Foundation (ISF) and
the US-Israel Binational Science Foundation (BSF) is gratefully
acknowledged.

\texttt{ \\Alex Lubotzky, \\ Einstein Institute of Mathematics, \\ The Hebrew University of Jerusalem, \\ Jerusalem 91904, \\Israel.\\e-mail: alexlub@math.huji.ac.il  \bigskip \\
Nikolay Nikolov, \\ School of Mathematics, \\ Tata Institute for Fundamental Research, \\ Mumbai 400 005, \\ India. \\e-mail: nikolov@math.tifr.res.in }

\begin{thebibliography}{99}
\bibitem{bt} A. Borel, J. Tits, El\'{e}ments unipotents et sousgroupes paraboliques des groupes r\'{e}ductifs I., Inv. Math. 12 (1971), 97-104.

\bibitem{BGLM} M. Burger, T. Gelander, A. Lubotzky, S. Mozes, Counting hyperbolic manifolds, Geometric and Functional Analysis (GAFA), Vol. 12 (2002) 1161-1173.
\bibitem{cooperstein}
B. N. Cooperstein, Minimal degree for a permutation representation of a classical group.
Israel J. Math. 30 (1978), no. 3, 213--235.

\bibitem{Edhan} O. Edhan. Counting congruence subgroups, M.Sc. Thesis, Hebrew University of Jerusalem, 2003.

\bibitem{GLP} D. Goldfeld, A. Lubotzky, L. Pyber, Counting congruence subgroups,
\emph{this journal.}

\bibitem{GLNP} D. Goldfeld, A. Lubotzky, N. Nikolov, L. Pyber, Counting primes, groups and manifolds, \emph{to appear in Proc. Nat. Acad. Sci. U.S.A.} 

\bibitem{cfsg} D. Gorenstein, R. Lyons, R. Solomon,
The Classification of the Finite Simple Groups. Number 3:
Almost simple $K$-groups. Mathematical Surveys and Monographs,
AMS, Providence, 1998.

\bibitem{KL} P. Kleidman, M. Liebeck, The Subgroup Structure of the Finite Classical Groups. LMS Lecture Note Series, 129. CUP, Cambridge, 1990.

\bibitem{lpink} M. Larsen, R. Pink, Finite subgroups of algebraic groups, J. Amer.
Math. Soc. \emph{to appear.}

\bibitem{LP} M. Liebeck, L. Pyber. Finite linear groups and bounded generation,
Duke Math. J. 107 (2001), no. 1, 159-171.

\bibitem{lss} M. Liebeck, J. Saxl, G. Seitz, Subgroups of maximal rank in finite exceptional groups of Lie type, Proc. London Math. Soc.(3) 65 (1992), 297-325.


\bibitem{ls} M. Liebeck, A. Shalev, Fuchsian groups, coverings of Riemann surfaces, subgroup growth, random quotients and random walks, \emph{preprint}.


\bibitem{lu} A. Lubotzky, Free quotients and the first Betti number of some hyperbolic manifolds, Transformation groups 1 (1996) 71-82.

\bibitem{LS} A. Lubotzky, D. Segal, Subgroup Growth, Progr. Math. 212, Birkh\"{a}user, Boston, 2003.

\bibitem{LM} A. Lubotzky, T. Weigel, Lattices of minimal covolume in $\mathrm{SL}_2$ over local fields, Proc. Lond. Math. Soc. 78 (1999) 283-333.

\bibitem{margulis} G. Margulis, Discrete Subgroups of Semisimple Lie Groups, Ergebnisse der Math. 17, Springer-Verlag, 1991.

\bibitem{MP} T. W. M\"{u}ller, J.-C. Puchta, Character theory of symmetric groups, subgroup growth of Fuchsian groups and random walks, \emph{preprint}.

\bibitem{pr} V. Platonov, A. Rapinchuk, Algebraic Groups and Number Theory, Pure and Applied Mathematics 139, Academic Press, Boston, 1994.

\bibitem{serre} J-P. Serre, Le probl\`{e}me des groupes de congruence pour $\mathrm{SL}_2$, Ann. of Math. (2) 92, (1970), 489-527.
\end{thebibliography}
\end{document}